\documentclass[11pt]{amsart}
\usepackage{amsmath}
\usepackage{amssymb}
\usepackage{amsfonts}

\newcommand{\C}{\mathbb{C}}
\newcommand{\Z}{\mathbb{Z}}

\newcommand{\N}{\mathbb{N}}
\newcommand{\R}{\mathbb{R}}

\newcommand{\h}{{\mathcal H}}
\newcommand{\lh}{{\mathcal L}({\mathcal H})}
\newcommand{\ran}{{\rm ran}}
\newcommand{\dist}{{\rm dist}}

\DeclareMathSymbol{\subsetneqq}{\mathbin}{AMSb}{36}

\newtheorem{th1}{{\bf Theorem}}[section]
\newtheorem{thm}[th1]{{\bf Theorem}}
\newtheorem{lem}[th1]{{\bf Lemma}}
\newtheorem{prop}[th1]{{\bf Proposition}}
\newtheorem{cor}[th1]{{\bf Corollary}}

\theoremstyle{remark}
\newtheorem{rem}[th1]{Remark}

\theoremstyle{definition}

\newtheorem{exm}[th1]{Example}
\newtheorem{ques}{Question}

\begin{document}
\title{ON THE LOCAL SPECTRAL PROPERTIES OF WEIGHTED SHIFT OPERATORS}
\author{A. Bourhim}
\address{Abdus Salam ICTP,
Mathematics Section, 11 Miramare, 34100 Trieste, Italy }
\email{bourhim@ictp.trieste.it}
\thanks{The research author was supported in part by the Abdus Salam International
Centre for Theoretical Physics, Trieste, Italy. The author is
grateful to the staff of Mathematics section of ICTP and Professor
C. E. Chidume for their generous hospitality\\ The present revised 
version was resubmitted to Studia Mathematica on January 03, 2003}


\begin{abstract}
In this paper, we study the local spectral properties for both
unilateral and bilateral weighted shift operators.
\end{abstract}

\maketitle

\section{Introduction}
\setcounter{equation}{0}

The weighted shift operators have been proven to be a very
interesting collection of operators, providing examples and
counter-examples to illustrate many properties of operators. The
basic facts concerning their spectral theory have been fully
investigated, and can be found in A. Shields' excellent survey
article \cite{shi}. The description of their local spectra, and
the characterization of their local spectral properties is
natural, and there has been some progress in the literature. In K.
Laursen, and M. Neumann's remarkable monograph \cite{ln}, it is
shown that each local spectrum of a unilateral weighted shift
operator $S$ contains a closed disc centered at the origin.
Sufficient conditions are given for a unilateral weighted shift
operator $S$ to satisfy Dunford's condition $(C)$ and for $S$ not
to possess Bishop's property $(\beta)$. It is also shown that a
unilateral weighted shift operator is decomposable if and only if
it is quasi-nilpotent. But, it is unknown which unilateral or
bilateral weighted shift operators satisfy Dunford's condition
$(C)$ or possess Bishop's property $(\beta)$. It is also unknown
which bilateral weighted shift operators are decomposable. The
main goal of this paper is to study the local spectral properties
of both unilateral and bilateral weighted shift operators. We
study the question of determining the local spectra of a weighted
shift $S$ in terms of the weight sequence defining $S$. Several
necessary and sufficient conditions for $S$ to satisfy Dunford's
condition $(C)$ or Bishop's property $(\beta)$ are given.

A part of the material of this paper is contained in the Ph.D thesis of the author, \cite{phd},
that was defended on April 2001.
After the first version of this paper, \cite{bourhim}, has been submitted for publication and
circulated among some specialists in the domain, the author learned from Professor M. M.
Neumann that some results of \cite{bourhim} had also appeared in \cite{mmn}. However, the
scope of the present paper is wider and the arguments of the proofs are simpler that the ones
given in \cite{phd} and \cite{mmn}. In \cite{mmn}, the authors established some  results abscent
in the present paper about inner and outer radii of arbitrary Banach space operator and gave
growth condition for a Banach space operator  to possess Bishop's property $(\beta)$.

Throughout this paper, $\h$ shall denote a complex Hilbert space,
and $\lh$ the algebra of all linear bounded  operators on $\h$.
For an operator $T\in \lh$, let $T^*$, $\sigma(T)$, $\sigma_{ap}(T)$,
$\sigma_p(T)$, $W(T)$, $r(T)$, and $w(T)$ denote the adjoint,
the spectrum, the approximate point spectrum, the point spectrum,
the numerical range, the spectral radius, and the numerical radius,
respectively, of $T$. Let $m(T):=\inf\{\|Tx\|:\|x\|=1\}$ denote the lower
bound of $T$, and let $r_1(T)$ denote $\sup\limits_{n\geq 1}[m(T^n)]^\frac{1}{n}$
which equals $\lim\limits_{n\to +\infty}[m(T^n)]^\frac{1}{n}$. An
operator $T\in \lh$ is said to have the  {\it single-valued
extension property} provided that for every non-empty open set
$U\subset\C$, the equation $(T-\lambda) f(\lambda)=0$ admits the
zero function $f\equiv 0$ as a unique analytic solution on $U$.
The {\it local resolvent set} $\rho_{_T}(x)$ of $T$ at a point
$x\in \h$ is the union of all open subsets $U\subset \C$ for which
there exists an analytic function $f:U\to \h$ such that
$(T-\lambda)f(\lambda)=x$ on $U$. The complement in $\C$ of
$\rho_{_T}(x)$  is called the {\it local spectrum} of $T$ at $x$
and will be denoted by $\sigma_{_T}(x)$. Recall that the {\it local
spectral radius} of $T$ at $x$ is given by $r_{_T}(x):=\limsup\limits_{n\to+\infty}\|T^nx\|^{\frac{1}{n}}$ and equals
$\max\{|\lambda|:\lambda\in\sigma_{_T}(x)\}$ whenever $T$ has the single-valued extension property.
Moreover, if $T$ enjoys this property  then for every $x\in\h$, there exists a unique maximal analytic solution $\widetilde{x}(.)$ on
$\rho_{_T}(x)$ for which $(T-\lambda)\widetilde{x}(\lambda)=x$ for
all $\lambda\in\rho_{_T}(x)$, and satisfies
$\widetilde{x}(\lambda)=-\sum\limits_{n\geq0}\frac{T^nx}{\lambda^{n+1}}$
on $\{\lambda\in\C: |\lambda|>r_{_T}(x)\}$. If in addition $T$ is invertible,
then $\widetilde{x}(\lambda)=\sum\limits_{n\geq
1}\lambda^{n-1}T^{-n}x$ on the open disc
$\{\lambda\in\C:|\lambda|<\frac{1}{r_{_{T^{-1}}}(x)}\}$. For a
closed subset $F$ of $\C$, let
$\h_{_T}(F):=\{x\in\h:\sigma_{_T}(x)\subset F\}$ be the
corresponding analytic spectral subspace; it is a
$T-$hyperinvariant subspace, generally non-closed in $\h$. The
operator $T$ is said to satisfy {\it Dunford's condition }$(C)$ if
for every closed subset $F$ of $\C$, the linear subspace
$\h_{_T}(F)$ is closed. Let $U$ be an open subset of $\C$ and let
${\mathcal O}(U,\h)$ be the space of analytic $\h-$valued
functions on $U$. Equipped with the topology of uniform
convergence on compact subsets of $U$, the space ${\mathcal
O}(U,\h)$ is a Fr\'echet space. Note that every operator $T\in
\lh$ induces a continuous mapping $T_{_U}$ on ${\mathcal O}(U,
\h)$ defined by $T_{_U}f(\lambda)=(T-\lambda)f(\lambda)$ for every
$f\in {\mathcal O}(U,\h)$ and $\lambda \in U$. An operator $T\in \lh$ is said to possess {\it Bishop's property $(\beta)$} provided that, for every
open subset $U$ of $\C$, the mapping $T_{_U}$ is injective and has a closed
range. Equivalently, if for every open subset $U$ of
$\C$ and for every sequence $(f_n)_n$ of ${\mathcal O}(U, H)$, the convergence of $(T_{_U}f_n)_n$ to $0$ in ${\mathcal O}(U, \h)$ should always
entail the convergence to $0$ of the sequence $(f_n)_n$ in
${\mathcal O}(U, \h)$. It is known that hyponormal
operators, $M-$hyponormal operators and more generally subscalar
operators possess Bishop's property $(\beta)$. It is also known
that Bishop's property ($\beta$) implies Dunford's condition $(C)$
and it turns out that the single-valued extension property follows
from Dunford's condition $(C)$. An operator $T\in\lh$ is said to
be {\it decomposable} if, for every open cover of $\C$ by two open
subsets $U_1$ and $U_2$, there exist $T-$invariant closed linear
subspaces $\h_1$ and $\h_2$ such that $\h=\h_1+\h_2$ and
$\sigma(T{_{|\h_i}})\subset U_i$ for $i=1,~2.$ Immediate examples
of decomposable operators are provided by the normal operators and
more generally by the spectral operators in the sense of Dunford.
It is known that an operator $T\in\lh$ is decomposable if and only
if both $T$ and $T^*$ possess Bishop's property $(\beta)$. For
thorough presentations of the local spectral theory, we refer to
the monographs \cite{cf}, and \cite{ln}.

Throughout this paper, let $S$ be a weighted shift operator on
$\h$ with a positive bounded weight sequence $(\omega_n)_{n}$,
that is $$Se_n=\omega_ne_{n+1}~\forall n,$$ where $(e_n)_{n}$ is
an orthonormal basis of $\h$. If the index $n$ runs over the
non-negative integers (resp. all integers), then $S$ is called a
{\it unilateral} (resp. {\it  bilateral}) {\it weighted shift }.

Before proceeding, we would like to recall some simple but useful
remarks which will be repeatedly used in the sequel.
\begin{itemize}
\item[$(R_1)$] For every $\alpha\in\C$, $|\alpha|=1$, we have
$(\alpha S)U_\alpha=U_\alpha S,$ where $U_\alpha$ is the unitary
operator defined on $\h$ by $U_\alpha e_n=\alpha^ne_n,~\forall n$.
\item[$(R_2)$] Let $K$ be a non-empty compact subset of $\C$. If $K$ is connected and invariant by circular symmetry about the origin,
then there are two real numbers $a$, and $b$ with $0\leq a\leq b$
such that $K=\{\lambda\in\C:a\leq|\lambda|\leq b\}.$
\item[$(R_3)$] Let $T$ and $R$ in $\lh$ such that $R$ is invertible. Then for every $x\in\h$, we have
$\sigma_{_T}(x)=\sigma_{_{RTR^{-1}}}(Rx).$
\end{itemize}

Finally, we need to introduce the local functional calculus,
developed by L. R. Williams in \cite{wil}, of operators having the
single-valued extension property; this local functional calculus
extends the Riesz functional calculus for operators. P. McGuire
has also studied this local functional calculus, but only for
operators with no eigenvalues (see \cite{MC}). Let $U$ be an open
subset of $\C$ and $K$ be a compact subset of $U$. Let
$(\gamma_i)_{1\leq i\leq n}$ be a finite family of disjoint closed
rectifiable Jordan curves in $U\backslash K$. The formal sum,
$\gamma:=\gamma_1+...+\gamma_n$, of the curves is called {\it an
oriented envelope of} $K$ in $U$ if its winding number equals $1$
on $K$ and $0$ on $\C\backslash U$, i.e., $$
\eta_\gamma(\lambda):=\frac{1}{2\pi
i}\oint_\gamma\frac{1}{\alpha-\lambda}d\alpha=\left\{
\begin{array}{lll}
1&\mbox{ for every }\lambda\in K\\
\\
0&\mbox{ for every }\lambda\in \C\backslash U
\end{array}
\right. $$ Let $T\in\lh$, and let $v\in\h$. Assume that $T$ has
the single-valued extension property. For an analytic complex
valued function $f$ on a neighborhood $U$ of $\sigma_{_T}(v)$, we
define a vector $f(T,v)$ by the equation $$f(T,v):=\frac{-1}{2\pi
i}\oint_\gamma f(\lambda)\widetilde{v}(\lambda)d\lambda,$$ where
$\gamma$ is an oriented envelope of $\sigma_{_T}(v)$ in $U$. Note
 that $f(T,v)$ is independent of the choice of $\gamma$. We end
this introduction by quoting, without proof, the following theorem
from \cite{wil} which will be used in the sequel.
\begin{thm}\label{williams}
Let $T\in\lh$, and let $v\in\h$. Assume that $T$ has the
single-valued extension property. If $f$ is an analytic function
on a neighborhood $O$ of $\sigma_{_T}(v)$ and is not identically
zero on each connected component of $O$, then the following
statements hold
\begin{itemize}
\item[$(a)$]$Z_{_T}(f,v):=\{\lambda\in\sigma_{_T}(v):f(\lambda)=0\}$ is a finite set.
\item[$(b)$]$\sigma_{_T}(v)=\sigma_{_T}(f(T,v))\cup Z_{_T}(f,v).$
\end{itemize}
Moreover, if $\sigma_p(T)=\emptyset$, then
$\sigma_{_T}(v)=\sigma_{_T}(f(T,v))$.
\end{thm}

\section{Preliminaries}

In this section, we assemble some simple but useful results which will be needed in the sequel.
These results remain valid in the general setting of Banach space operators.

Let $\lambda_0\in\C$; recall that an operator $T\in \lh$ is said to
have {\it the single-valued extension property }
(resp. possess {\it Bishop's property $(\beta)$}) at $\lambda_0$
if there is an open neighbourhood $V$ of $\lambda_0$ such that for every open subset $U$ of $V$, the mapping $T_{_U}$ is injective (resp. is injective and has a closed range).
Obviously, $T$ has the single-valued extension property at
$\lambda_0$ whenever $T$ possesses Bishop's property $(\beta)$ at $\lambda_0$. Note that $T$ possesses Bishop's property $(\beta)$ at $\lambda_0$ if and only if there is an open neighbourhood $V$ of
$\lambda_0$ such that for every open subset $U$ of $V$ and for
every sequence $(f_n)_n$ of ${\mathcal O}(U, H)$, the convergence
of $(T_{_U}f_n)_n$ to $0$ in ${\mathcal O}(U, \h)$ should always
entail the convergence to $0$ of the sequence $(f_n)_n$ in
${\mathcal O}(U, \h)$. Moreover, if $T$ possesses Bishop's property $(\beta)$ at any point $\lambda\in\C$
then $T$ possesses Bishop's classical property $(\beta)$.

For an arbitrary $T\in\lh$, we shall denote
$$\sigma_{\beta}(T):=\big\{\lambda\in\C:T\mbox{ fails to possess Bishop's property }(\beta)\mbox{ at }\lambda\big\}.$$
It is clear that $\sigma_{\beta}(T)$ is a closed subset of $\C$ contained in $\sigma(T)$. In fact more can be established in the following result.

\begin{prop}\label{kim1}For any operator $T\in\lh$, we have
$\sigma_{\beta}(T)\subset\sigma_{ap}(T).$
\end{prop}
\begin{proof}Let $\lambda_0\not\in\sigma_{ap}(T)$, we shall show that $T$ possesses Bishop's
property $(\beta)$ at $\lambda_0$. By replacing $T$ with $T-\lambda_0$, we may assume without loss
of generality that $\lambda_0=0$. In this case, there is $R\in\lh$ such that
$RT=1$. Therefore, for every $\lambda\in\C$ and for every $x\in\h$, we have
\begin{eqnarray*}
\|(T-\lambda)x\|&=&\|(1-\lambda R)Tx\|\\
&\geq&m(1-\lambda R)m(T)\|x\|\\
&=&\frac{1}{\|R\|}m(1-\lambda R)\|x\|.
\end{eqnarray*}
If we set $V:=\{\lambda\in\C:|\lambda|<\frac{1}{\|R\|}\}$ then for every
$\lambda\in V$, $1-\lambda R$ is invertible. Hence,
for every $\lambda\in V$ and for every $x\in\h$, we have
\begin{equation}\label{kim3}
\|x\|\leq\|R\|\|(1-\lambda R)^{-1}\|\|(T-\lambda)x\|.
\end{equation}
Now, let $U$ be an open subset of the open disc $V$ and let
$(f_n)_{n\geq 0}$ be a sequence of ${\mathcal O}(U,\h)$ such that
$(T_{_U}f_n)_n$ converges to $0$ in ${\mathcal O}(U,\h)$. Let $K$
be a compact subset of $U$, it follows from (\ref{kim3}) that
$$\|f_n(\lambda)\|\leq\|R\|\max\big\{\|(1-\mu R)^{-1}\|:\mu\in K\big\}\|(T-\lambda)f_n(\lambda)\|,~~\forall \lambda\in K.$$
This shows that the sequence $(f_n)_n$ converges to $0$ uniformly
on $K$; and so, $S$ possesses Bishop's property $(\beta)$ any point of $V$. In particular,
$S$ possesses Bishop's property $(\beta)$ at $\lambda_0=0$. This finishes the proof.
\end{proof}

For every operator $T\in\lh$, we shall denote
$$\Re(T):=\big\{\lambda\in\C:T\mbox{ does not have the single-valued extension property at }\lambda\big\}.$$
It is an open subset of $\C$ contained in $\sigma(T)$. The description of
$\Re(T)$ for several operators $T$ can be found in \cite{amn} and \cite{abdou}. Note that
$$\overline{\sigma(T)\backslash\sigma_{ap}(T)}\subset\Re(T^*)$$
for every operator $T\in\lh$. This containment fails to be equality in general as counter-example can
be given by a unilateral weighted shift operator. In \cite{abdou}, it is shown that if $T$ is a cyclic operator
possessing Bishop's property $(\beta)$ then $$\Re(T^*)=\overline{\sigma(T)\backslash\sigma_{ap}(T)}.$$
In fact, this result remain valid for every rationally cyclic injective operator possessing Bishop's property $(\beta)$.
Recall that an operator (resp. injective operator) $T\in\lh$ is said to be  {\it cyclic } (resp. {\it rationally cyclic })
if there is a vector $x\in\h$ (resp. $x\in T^\infty\h:=\bigcap\limits_{n\geq0}T^n\h$) such that the closed linear subspace
generated by $\{T^nx:n\geq0\}$ (resp. $\{T^nx:n\in\Z\}$) is dense in $\h$.

The proof of the following result is an adaptation of the one of theorem 3.1 of \cite{abdou} and will be omitted.

\begin{prop}\label{rationally}
Assume that $T\in\lh$ is a rationally cyclic injective operator. If $T$ possesses Bishop's property $(\beta)$, then
$\Re(T^*)=\overline{\sigma(T)\backslash\sigma_{ap}(T)}.$
\end{prop}

The following two lemmas are easy to verify, we omit the proof of the next one.

\begin{lem}\label{kim2}Let $\h_1$ and $\h_2$ be complex Hilbert spaces, and let
$$M_C=\left[
\begin{array}{cc}
A&C\\
0&B
\end{array}
\right]$$
be a $2\times2$ upper triangular operator matrix on $\h_1\oplus\h_2$. If both $A$ and
$B$ possess Bishop's property $(\beta)$ at a point $\lambda_0\in\C$, then $M_C$ possesses
Bishop's property $(\beta)$ at $\lambda_0$.
\end{lem}

\begin{lem}\label{kim33}Let $\h_1$ and $\h_2$ be complex Hilbert spaces, and let
$$M_C=\left[
\begin{array}{cc}
A&C\\
0&B
\end{array}
\right]$$
be a $2\times2$ upper triangular operator matrix on $\h_1\oplus\h_2$. If $B$ have the single-valued extension
property, then for every $x\in\h_1$, we have
$$\sigma_{_{M_C}}(x\oplus0)=\sigma_{_A}(x).$$
\end{lem}
\begin{proof}Let $x\in\h_1$.

Let $U$ be an open subset of $\C$ and $f:U\to\h_1$ be an analytic function such that
$$(A-\lambda)f(\lambda)=x,\mbox{ for every }\lambda\in U.$$
If for every $\lambda\in U$, we set $F(\lambda)=f(\lambda)\oplus0,$
then
$$(M_C-\lambda)F(\lambda)=x\oplus0,\mbox{ for every }\lambda\in U.$$
This shows that
$$\sigma_{_{M_C}}(x\oplus0)\subset\sigma_{_A}(x).$$

Conversely, let $U$ be an open subset of $\C$ and $F:U\to\h$ be an analytic function such that
$$(M_C-\lambda)F(\lambda)=x\oplus0,\mbox{ for every }\lambda\in U.$$
Write $F=f_1\oplus f_2$, where $f_1:U\to\h_1$ and $f_2:U\to\h_2$ are analytic functions.
For every $\lambda\in U$, we have
$$(A-\lambda)f_1(\lambda)+Cf_2(\lambda)=x\mbox{ and }(B-\lambda)f_2(\lambda)=0.$$
Since $B$ has the single-valued extension property, we have $f_2\equiv0$.
Hence,
$$(A-\lambda)f_1(\lambda)=x\mbox{ for every }\lambda\in U.$$
This shows that $$\sigma_{_A}(x)\subset\sigma_{_{M_C}}(x\oplus0).$$
And the proof is complete.
\end{proof}

Obviously, it is seen from $(R_1)$ that $\sigma(S)$, $\sigma_{ap}(S)$ and $\sigma_p(S)$ have circular
symmetry about the origin. However, the fact that $\sigma_\beta(S)$ is circularly symmetric deserves a proof.

\begin{lem}\label{kim5}
If $S$ is a unilateral or bilateral weighted shift on $\h$, then $\sigma_\beta(S)$ is
invariant under circular symmetry about the origin.
\end{lem}
\begin{proof}
Assume that $S$ possesses Bishop's property
$(\beta)$ at a point $\lambda_0\in\C$. There is $\delta>0$ such that for every open subset $U$ of the disc
$B:=\{\lambda\in\C:|\lambda-\lambda_0|<\delta\}$, and for every sequence
$(f_n)_n$ of ${\mathcal O}(U,\h)$, the convergence of $(S_Uf_n)_n$
to $0$ in ${\mathcal O}(U,\h)$ implies the convergence of
$(f_n)_n$ to $0$ in ${\mathcal O}(U,\h)$.
Let $\alpha\in\C$, $|\alpha|=1$; it should be seen that $S$ possesses Bishop's property $(\beta)$ at $\alpha\lambda_0$.
Indeed, let $V$ be an open subset of $\alpha B$, and let $(g_n)_n$ be a sequence of ${\mathcal O}(V,\h)$ for which $(S_Vg_n)_n$ converges to $0$ in ${\mathcal O}(V,\h)$. Let $K$ be a compact subset of $V$; we shall show that
$\lim\limits_{n\to+\infty}\big[\sup\limits_{\lambda\in
K}\|g_n(\lambda)\|\big]=0$. Note that in view of $(R_1)$, we have
\begin{equation}\label{kkim}
\|(S-\alpha\lambda)x\|=\|(S-\lambda)U_\alpha
x\|,~\forall\lambda\in\C,~\forall x\in\h.
\end{equation}
And so, if we set
$f_n(\lambda):=g_n(\alpha\lambda)$,~$\lambda\in\alpha^{-1}V$, then
it follows from (\ref{kkim}) that
\begin{eqnarray*}
\sup\limits_{\lambda\in \alpha^{-1}K}\|(S-\lambda)U_\alpha
f_n(\lambda)\|&=&\sup\limits_{\lambda\in
\alpha^{-1}K}\|(S-\alpha\lambda)f_n(\lambda)\| \\
&=&\sup\limits_{\lambda\in
\alpha^{-1}K}\|(S-\alpha\lambda)g_n(\alpha\lambda)\|\\
&=&\sup\limits_{\lambda\in K}\|(S-\lambda)g_n(\lambda)\|.
\end{eqnarray*}
As $\lim\limits_{n\to+\infty}\big[ \sup\limits_{\lambda\in
K}\|(S-\lambda)g_n(\lambda)\|\big]=0$, and $S$ possesses Bishop's
property $(\beta)$ at $\lambda_0$, we have
$$\lim\limits_{n\to+\infty}\big[\sup\limits_{\lambda\in
\alpha^{-1}K}\|U_\alpha f_n(\lambda)\|\big]=0.$$ Since $U_\alpha$
is an unitary operator, we have
$$\lim\limits_{n\to+\infty}\big[\sup\limits_{\lambda\in
K}\|g_n(\lambda)\|\big]=\lim\limits_{n\to+\infty}\big[\sup\limits_{\lambda
\in\alpha^{-1}K}\|f_n(\lambda)\|\big]=0. $$ And the proof is
complete.
\end{proof}

\section{Local spectral properties of unilateral weighted shift operators}

In dealing with unilateral weighted shift operators, let
$(\beta_n)_{n\geq0}$ be the sequence given by $$ \beta_n=\left\{
\begin{array}{lll}
\omega_0...\omega_{n-1}&\mbox{if }n>0\\
\\
1&\mbox{if }n=0
\end{array}
\right. $$ and let us associate to $S$ the following quantities,
$$r_2(S):=\liminf\limits_{n\to\infty}\big{[}\beta_n
\big{]}^{\frac{1}{n}}, \mbox{   and   }
r_3(S):=\limsup\limits_{n\to\infty}\big{[}  \beta_n
\big{]}^{\frac{1}{n}}.$$ Note that,
$$r(S)=\lim\limits_{n\to\infty}\bigg{[}\sup\limits_{k\geq 0}
\frac{\beta_{n+k}}{\beta_k}  \bigg{]}^{\frac{1}{n}},~
r_1(S)=\lim\limits_{n\to\infty}\bigg{[}\inf\limits_{k\geq 0}
\frac{\beta_{n+k}}{\beta_k}\bigg{]}^{\frac{1}{n}},$$ and
$$m(S)\leq r_1(S)\leq r_2(S)\leq r_3(S)\leq r(S)\leq
w(S)\leq\|S\|.$$

\subsection{On the local spectra of unilateral weighted shift operators}

The following result can be deduced from proposition 1.6.9 of
\cite{ln}; we give here a direct proof.

\begin{prop}\label{r_2}
For every non-zero $x\in\h$, we have
$$\{\lambda\in\C:|\lambda|\leq r_2(S)\}\subset\sigma_{_S}(x).$$
\end{prop}
\begin{proof}As $\bigcap\limits_{n\geq 0}S^n\h=\{0\}$, we get that $0\in\sigma_{_S}(x)$
for every $x\in\h\backslash\{0\}$. Thus, we may assume that
$r_2(S)>0$. Let $O:=\{\lambda\in\C:|\lambda|<r_2(S)\}$, and
consider the following analytic $\h$--valued function on $O$,
$$k(\lambda):=\sum\limits_{n=0}^{+\infty}\frac{\lambda^n}{\beta_n}e_n.$$
It is easy to see that $(S-\lambda)^*k(\overline{\lambda})=0$ for
every $\lambda\in O$. Now, let $x=\sum\limits_{n\geq0}a_ne_n\in\h$
such that $O\cap\rho_{_S}(x)\not=\emptyset$. So, for every
$\lambda\in O\cap\rho_{_S}(x),$ we have
\begin{eqnarray*}
\widehat{x}(\lambda)&:=&\langle x;k(\overline{\lambda})\rangle\\
&=&\langle
(S-\lambda)\widetilde{x}(\lambda);k(\overline{\lambda})\rangle\\
&=&\langle
\widetilde{x}(\lambda);(S-\lambda)^*k(\overline{\lambda})\rangle\\
&=&0.
\end{eqnarray*}
Since
$\widehat{x}(\lambda)=\sum\limits_{n=0}^{+\infty}\frac{a_n}{\beta_n}\lambda^ne_n$
for every $\lambda\in O$, we see that $a_n=0$ for every $n\geq0$,
and therefore, $x=0$. Thus, the proof is complete.
\end{proof}

\begin{cor}\label{con}
For every $x\in\h$, $\sigma_{_S}(x)$ is a connected set.
\end{cor}
\begin{proof}Suppose for the sake of contradiction that there is a non-zero element $x$ of $\h$ such that $\sigma_{_S}(x)$
is disconnected. So, there is two non-empty disjoint compact
subsets $\sigma_1$ and $\sigma_2$ of $\C$ such that
$\sigma_{_S}(x)=\sigma_1\cup\sigma_2$. It follows from proposition
1.2.16 (g) of \cite{ln} that
$\h_{_S}\big(\sigma_{_S}(x)\big)=\h_{_S}(\sigma_1)\oplus\h_{_S}(\sigma_2)$.
Therefore, there are two non-zero elements $x_1$ and $x_2$ of $\h$
such that $x=x_1+x_2$ and $\sigma_{_S}(x_i)\subset \sigma_i$,
$i=1,$ $2$. In particular,
$\sigma_{_S}(x_1)\cap\sigma_{_S}(x_2)=\emptyset$, which is a
contradiction to proposition \ref{r_2}.
\end{proof}

We give here a simple proof of theorem 4 of \cite{shi} from the
point of view of the local spectral theory.

\begin{cor}\label{spec}The spectrum of $S$ is the disc $\{\lambda\in\C:|\lambda|\leq r(S)\}.$
\end{cor}
\begin{proof}By proposition 1.3.2 of
\cite{ln}, we have
$\sigma(S)=\bigcup\limits_{x\in\h\backslash\{0\}} \sigma_{_S}(x)$.
As each $\sigma_{_S}(x)$ is a connected set, and
$\{\lambda\in\C:|\lambda|\leq
r_2(S)\}\subset\bigcap\limits_{x\in\h\backslash\{0\}}\sigma_{_S}(x)\not=\emptyset$
(see proposition \ref{r_2}), we deduce that $\sigma(S)$ is also a
connected set containing the disc $\{\lambda\in\C:|\lambda|\leq
r_2(S)\}.$ On the other hand, from $(R_1)$ we see that $\sigma(S)$
has circular symmetry about the origin. So, by $(R_2)$,
$\sigma(S)$ must be the disc $\{\lambda\in\C:|\lambda|\leq
r(S)\}.$
\end{proof}
\begin{rem}\label{r_3}
Let $x=\sum\limits_{n\geq 0}a_ne_n$ be a non-zero element of $\h$.
It follows from proposition \ref{r_2}, that $r_2(S)\leq
r_{_S}(x)$. In fact, more can be shown:
\begin{equation}\label{rsx}
r_3(S)\leq r_{_S}(x)\leq r(S).
\end{equation}
Indeed, for every integer $n\geq0$, we have
\begin{eqnarray*}
S^nx&=&\sum\limits_{k\geq0}a_kS^ne_k\\
&=&\sum\limits_{k\geq0}a_k\frac{\beta_{n+k}}{\beta_k}e_{n+k}.
\end{eqnarray*}
Since there is an integer $k_0\geq0$ such that $a_{k_0}\not=0$, it
then follows that
$$|a_{k_0}|\frac{\beta_{n+k_0}}{\beta_{k_0}}\leq\bigg[\sum\limits_{k\geq0}|a_k|^2\bigg{|}
\frac{\beta_{n+k}}{\beta_k}\bigg{|}^2
\bigg]^{\frac{1}{2}}=\|S^nx\|,~\forall n\geq0.$$ Now, by taking
the $n$th root and then $\limsup$ as $n\to+\infty$, we get
$r_3(S)\leq r_{_S}(x)$.
\end{rem}

Let $x$ be a non-zero element of $\h$; from corollary \ref{con},
we see that every circle of radius $r$, $0\leq r\leq r_{_S}(x)$,
intersects $\sigma_{_S}(x)$. So, in view of remark \ref{r_3},
$\sigma_{_S}(x)$ may contain points which are not in the disc
$\{\lambda\in\C:|\lambda|\leq r_2(S)\}$. The next result gives a
complete description of local spectrum of $S$ at most of the
points in $\h$, and refines the local spectral inclusion in
proposition \ref{r_2}. For every non-zero $x=\sum\limits_{n\geq
0}a_ne_n\in\h$, we set
$$R_\omega(x):=\liminf\limits_{n\to+\infty}\big{|}
\frac{\beta_n}{a_n} \big{|}^{\frac{1}{n}}.$$ Note that,
$r_2(S)\leq R_\omega(x)\leq +\infty$ for every non-zero $x\in\h$.

\begin{thm}\label{lsp}For every non-zero $x=\sum\limits_{n\geq 0}a_ne_n\in\h$, we have:
\begin{itemize}
\item[$(a)$]If $r_3(S)<R_\omega(x)$, then $\sigma_{_S}(x)=\{\lambda\in\C:|\lambda|\leq r_3(S)\}.$
\item[$(b)$]If $r_3(S)\geq R_\omega(x)$, then $\{\lambda\in\C:|\lambda|\leq R_\omega(x)\}\subset\sigma_{_S}(x).$
\end{itemize}
\end{thm}
\begin{proof}Let us first show that
\begin{equation}\label{sle0}
\sigma_{_S}(e_0)=\{\lambda\in\C:|\lambda|\leq r_3(S)\}.
\end{equation}
Let $\alpha\in\C$, $|\alpha|=1$; by $(R_1)$, we have $\alpha
S=U_\alpha SU_\alpha^*$. And so, by $(R_3)$, we have
$\sigma_{_S}(e_0)=\sigma_{_{\alpha S}}(U_\alpha e_0)$; hence,
$\sigma_{_S}(e_0)=\alpha\sigma_{_S}(e_0)$. This shows that
$\sigma_{_S}(e_0)$ has circular symmetry about the origin. As
$\sigma_{_S}(e_0)$ is a connected set containing the disc
$\{\lambda\in\C:|\lambda|\leq r_2(S)\}$, it follows from $(R_2)$
that $\sigma_{_S}(e_0)$ must be the disc
$\{\lambda\in\C:|\lambda|\leq r_{_S}(e_0)\}$. Therefore,
(\ref{sle0}) holds, since $r_{_S}(e_0)=r_3(S)$.

Next, keep in mind that for every
$\lambda\in\rho_{_S}(e_0)=\{\lambda\in\C:|\lambda|>r_3(S)\}$, we
have
\begin{eqnarray*}
\widetilde{e_0}(\lambda)&=&-\sum\limits_{n\geq0}\frac{S^ne_0}{\lambda^{n+1}}\\
&=&-\sum\limits_{n\geq0}\frac{\beta_n}{\lambda^{n+1}}e_n.
\end{eqnarray*}
Now, we are able to prove the first statement.

$(a)$ Suppose that $r_3(S)<R_\omega(x)$. Note that the following
function $f(\lambda):=\sum\limits_{n\geq
0}\frac{a_n}{\beta_n}\lambda^n$ is analytic on the neighborhood,
$\{\lambda\in\C:|\lambda|<R_\omega(x)\}$, of $\sigma_{_S}(e_0)$.
Let $r$ be a real number such that $r_3(S)<r<R_\omega(x)$, we have
\begin{eqnarray*}
f(S,e_0)&=&\frac{-1}{2\pi
i}\oint_{|\lambda|=r}f(\lambda)\widetilde{e_0}(\lambda)d\lambda\\
&=&\frac{-1}{2\pi
i}\oint_{|\lambda|=r}f(\lambda)\bigg[-\sum\limits_{n\geq
0}\frac{\beta_n}{\lambda^{n+1}}e_n\bigg]d\lambda\\
&=&\sum\limits_{n\geq 0}\bigg[\frac{1}{2\pi
i}\oint_{|\lambda|=r}\frac{f(\lambda)}{\lambda^{n+1}} d\lambda
\bigg]\beta_ne_n\\ &=&x.
\end{eqnarray*}
And so, by theorem \ref{williams}, we have
$$\sigma_{_S}(x)=\sigma_{_S}(e_0)=\{\lambda\in\C:|\lambda|\leq
r_3(S)\}.$$

$(b)$ Suppose that $r_3(S)\geq R_\omega(x)$. If
$r_2(S)=R_\omega(x)$, then according to proposition \ref{r_2}
there is nothing to prove; thus we may suppose that
$r_2(S)<R_\omega(x)$. For each $n\geq0$, let
$$A_n(\lambda)=-\frac{a_0}{\lambda^{n+1}}\beta_n-\frac{a_1}{\lambda^{n}}\frac{\beta_n}{\beta_1}-
\frac{a_2}{\lambda^{n-1}}\frac{\beta_n}{\beta_2}-...-\frac{a_n}{\lambda},~\lambda\in\C\backslash\{0\},$$
and
$$B_n(\lambda)=a_0+\frac{a_1}{\beta_1}\lambda+\frac{a_2}{\beta_2}\lambda^2+...
+\frac{a_n}{\beta_n}\lambda^n,~\lambda\in\C.$$ We have,
\begin{equation}\label{AB}
A_n(\lambda)=-\frac{\beta_n}{\lambda^{n+1}}B_n(\lambda),~\lambda\in\C\backslash\{0\}.
\end{equation}
By writing $\widetilde{x}(\lambda):=\sum\limits_{n\geq
0}\widetilde{A}_n(\lambda)e_n,~\lambda\in\rho_{_S}(x)$, we get
from the following equation,
$$(S-\lambda)\widetilde{x}(\lambda)=x,~\lambda\in\rho_{_S}(x),$$
that for every $\lambda\in\rho_{_S}(x)$, we have $$ \left\{
\begin{array}{lll}
-\lambda\widetilde{A}_0(\lambda)=a_0\\
\\
\widetilde{A}_n(\lambda)\omega_n-\lambda
\widetilde{A}_{n+1}(\lambda)=a_{n+1}&\mbox{ for every }n\geq 0.
\end{array}
\right. $$ Therefore, for every $n\geq0$, we have
$$\widetilde{A}_n(\lambda)=A_n(\lambda)=-\frac{a_0}{\lambda^{n+1}}\beta_n-\frac{a_1}{\lambda^{n}}\frac{\beta_n}{\beta_1}-
\frac{a_2}{\lambda^{n-1}}\frac{\beta_n}{\beta_2}-...-\frac{a_n}{\lambda},~\lambda\in\rho_{_S}(x).$$
Since
$\|\widetilde{x}(\lambda)\|^2=\sum\limits_{n\geq0}|\widetilde{A}_n(\lambda)|^2=\sum\limits_{n\geq0}|A_n(\lambda)|^2<+\infty$
for every $\lambda\in\rho_{_S}(x)$, it then follows that
\begin{equation}\label{lim}
\lim\limits_{n\to+\infty}A_n(\lambda)=0\mbox{ for every
}\lambda\in\rho_{_S}(x).
\end{equation}
We shall show that (\ref{lim}) is not satisfied for most of the
points in the open disc
$V(x):=\{\lambda\in\C:|\lambda|<R_\omega(x)\}$. It is clear that
the sequence $(B_n)_{n\geq0}$ converges uniformly on compact
subsets of $V(x)$ to the non-zero power series
$B(\lambda)=\sum\limits_{n\geq0}\frac{a_n}{\beta_n}\lambda^n$.
Now, let $\lambda_0\in V(x)\backslash\{0\}$ such that
$B(\lambda_0)\not=0$; there is $\epsilon>0$ and an integer $n_0$
such that $\epsilon<|B_n(\lambda_0)|$ for every $n\geq n_0$. On
the other hand, $|\lambda_0|<r_3(S)$, then there is a subsequence
$(n_k)_{k\geq0}$ of integers greater than $n_0$ such that
$|\lambda_0|^{n_k}<\beta_{n_k}$. Thus, (\ref{AB}) gives
$$|A_{n_k}(\lambda_0)|=|-\frac{\beta_{n_k}}{\lambda_0^{n_k+1}}B_{n_k}(\lambda_0)|\geq
\frac{\epsilon}{|\lambda_0|},\mbox{ for every }k\geq0.$$ And so,
by (\ref{lim}), $\lambda_0\not\in\rho_{_S}(x)$. Since the set of
zeros of $B$ is at most countable, we have
$\{\lambda\in\C:|\lambda|\leq R_\omega(x)\}\subset\sigma_{_S}(x).$
\end{proof}

\begin{rem}
Let $\h_0$ denote all the finite linear combinations of the
vectors $e_n,~(n\geq0)$; it is clearly a dense subspace of $\h$.
It follows from theorem \ref{lsp}$(a)$ that for every non-zero
$x\in\h_0$, we have $$\sigma_{_S}(x)=\{\lambda\in\C:|\lambda|\leq
r_3(S)\}.$$ This can also be deduced from theorem \ref{williams},
and (\ref{sle0}), since for every non-zero $x\in\h_0$, there is a
non-zero polynomial $p$ such that $x=p(S)e_0$.
\end{rem}

\begin{ques}
Does $\sigma_{_S}(x)=\{\lambda\in\C:|\lambda|\leq r_{_S}(x)\}$ for
every non-zero $x\in\h$?
\end{ques}
In view of proposition \ref{r_2}, and (\ref{rsx}), an interesting
special case of this question is suggested.
\begin{ques}\label{reso}Is $\{\lambda\in\C:|\lambda|\leq
r_3(S)\}\subset\sigma_{_S}(x)$ for every non-zero $x\in\h$?
\end{ques}

\subsection{Which unilateral weighted shift operators satisfy Dunford's condition $(C)$?}

The following result gives a necessary condition for the
unilateral weighted shift operator $S$ to enjoy Dunford's
condition $(C)$.

\begin{thm}\label{lah}If $S$ satisfies Dunford's condition $(C)$, then $r(S)=r_3(S)$.
\end{thm}
\begin{proof}
Suppose that $S$ satisfies Dunford's condition $(C)$, and let
$F:=\{\lambda\in\C:|\lambda|\leq r_3(S)\}$. It follows from the
statement $(a)$ of theorem \ref{lsp} that $\h_{_S}(F)$ contains a
dense subset of $\h$. As the subspace $\h_{_S}(F)$ is closed, we
get $\h_{_S}(F)=\h$; this means that, $\sigma_{_S}(x)\subset F$
for every $x\in\h$. And so,
$\sigma(S)=\bigcup\limits_{x\in\h}\sigma_{_S}(x)\subset F$ (see
proposition 1.3.2 of \cite{ln}). As $r_3(S)\leq r(S)$, and
$\sigma(S)=\{\lambda\in\C:|\lambda|\leq r(S)\}$, we get
$r(S)=r_3(S)$.
\end{proof}

An interesting special case occurs when the sequence
$([\beta_n]^{\frac{1}{n}})_{n\geq1}$ converges. Then by combining
proposition \ref{r_2}, and theorem \ref{lah}, we see that the
following statements are equivalent.
\begin{itemize}
\item[$(a)$] $\sigma_{_S}(x)=\sigma(S)$ for every non-zero $x\in\h$.
\item[$(b)$] $S$ satisfies Dunford's condition $(C)$.
\item[$(c)$] $r(S)=\lim\limits_{n\to+\infty}[\beta_n]^{\frac{1}{n}}$.
\end{itemize}
In considering the general case, we note that if $r(S)=r_3(S)$,
then $\sigma_{_S}(x)=\sigma(S)$ for most of the points $x$ in
$\h$. Thus, we conjecture that the following statements are
equivalent.
\begin{itemize}
\item[$(a)$] $\sigma_{_S}(x)=\sigma(S)$ for every non-zero $x\in\h$.
\item[$(b)$] $S$ satisfies Dunford's condition $(C)$.
\item[$(c)$] $r(S)=r_3(S)$.
\end{itemize}
Note that a positive answer to question \ref{reso} will prove this
conjecture.

\subsection{Which unilateral weighted shift operators possess Bishop's property $(\beta)$?}

In \cite{mm}, T. L. Miller and V. G. Miller have shown that if
$r_1(S)<r_2(S)$, then $S$ does not possess Bishop's property
$(\beta)$. In fact, this result could be obtained by combining
theorem 3.1 of \cite{yan} and theorem 10(ii) of \cite{shi}. But
the proof given in \cite{mm} is simple, and elegant. Here, we
refine this result as follows.

\begin{thm}\label{beta1}If $S$ possesses Bishop's property $(\beta)$,
then $r_1(S)=r(S)$. Conversely, if $r_1(S)=r(S)$, then either $S$
possesses Bishop's property $(\beta)$, or
$\sigma_{\beta}(S)=\{\lambda\in\C:|\lambda|=r(S)\}.$
\end{thm}
\begin{proof}
Suppose that $S$ possesses Bishop's property $(\beta)$. According
to the above discussion, we have
\begin{equation}\label{e1}
r_1(S)=r_2(S).
\end{equation}
On the other hand, it follows from \cite{AA} that $S$ is
power-regular, that is the sequence
$(\|S^nx\|^{\frac{1}{n}})_{n\geq0}$ converges for all $x\in\h$. In
particular, $(\|S^ne_0\|^{\frac{1}{n}})_{n\geq0}$ converges; so,
\begin{equation}\label{e2}
r_2(S)=r_3(S).
\end{equation}
As $S$ satisfies Dunford's condition $(C)$, we get from theorem
\ref{lah} that
\begin{equation}\label{e3}
r_3(S)=r(S).
\end{equation}
Therefore, the result follows from (\ref{e1}), (\ref{e2}), and
(\ref{e3}).

Conversely, suppose that $r_1(S)=r(S)$. Since $\sigma_{ap}(S)=\{\lambda\in\C:|\lambda|=r(S)\}$ (see theorem 1 of \cite{rid} or theorem 6 of \cite{shi}), the desired result
holds by combining proposition \ref{kim1} and lemma \ref{kim5}.
\end{proof}

In view of theorem \ref{beta1}, the following question is
suggested.

\begin{ques}\label{ques2}
Suppose that $r_1(S)=r(S)$. Does $S$ possess Bishop's property
$(\beta)$?
\end{ques}

The next two propositions give a condition on the weight
$(\omega_n)_{n\geq0}$ for which the corresponding weighted shift
possesses Bishop's property $(\beta)$.

\begin{prop}\label{w}If $m(S)=w(S)$, then $S$ possesses Bishop's property $(\beta)$.
\end{prop}
\begin{proof}It is well known that for every operator $T\in\lh$, and for every $\lambda$ not in the closure of $W(T)$, the norm of the
resolvent $(T-\lambda)^{-1}$ admits the estimate
\begin{equation}\label{w(T)}
\|(T-\lambda)^{-1}\|\leq\frac{1}{\dist(\lambda,W(T))}.
\end{equation}
If $m(S)=w(S)=0$, then $S$ is quasi-nilpotent; and so, by
proposition 1.6.14 of \cite{ln}, $S$ is decomposable. Therefore,
$S$ possesses Bishop's property $(\beta)$. Now, without loss of
generality suppose that $m(S)=w(S)=1$. It follows from proposition
16 of \cite{shi} that $$\{\lambda\in\C:|\lambda|<1\}\subset
W(S)\subset \{\lambda\in\C:|\lambda|\leq1\};$$ and so, by
(\ref{w(T)}) we have
$$\|(S-\lambda)^{-1}\|\leq\frac{1}{|1-|\lambda||},\mbox{ for every
} \lambda\in\C,~|\lambda|>1.$$ This shows that, for every
$x\in\h$, we have $$|1-|\lambda||\|x\|\leq\|(S-\lambda)x\|,\mbox{
for all } \lambda\in\C,~|\lambda|>1.$$ On the other hand, for
every $x\in\h$ and $\lambda\in\C$, $|\lambda|<1$, we have
\begin{eqnarray*}
\|(S-\lambda)x\|&\geq&\|Sx\|-|\lambda|\|x\|\\
&\geq&m(S)\|x\|-|\lambda|\|x\|\\ &\geq&(1-|\lambda|)\|x\|.
\end{eqnarray*}
Therefore, for every $\lambda\in\C\backslash\{\mu\in\C:|\mu|=1\}$
and for every $x\in\h$, we have $$|1-|\lambda||\|x\|\leq
\|(S-\lambda)x\|.$$ And so, it follows from theorem 1.7.1 of
\cite{ln} that $S$ possesses Bishop's property $(\beta)$.
\end{proof}

The next example gives a collection of unilateral weighted shift
operators possessing Bishop's property $(\beta)$ which show that
the assumption $m(S)=w(S)$ in proposition \ref{w} is not a
necessary condition.

\begin{exm}Let $a$ be a positive real number, and let $S_a$ be the unilateral weighted shift with the corresponding weight
$(\omega_n(a))_{n\geq0}$ given by $$\omega_n(a)=\left\{
\begin{array}{lll}
a&\mbox{if }n=0\\
\\
1&\mbox{if }n>0
\end{array}
\right.$$ It is clear that $m(S_a)=1$, and $S_a$ possesses
Bishop's property $(\beta)$ since $S_a$ and the unweighted shift
operator are unitarily equivalent. On the other hand, it was shown
in \cite{BS} (see also proposition 2 of \cite{tin}) that
\begin{itemize}
\item[$(a)$]If $a\leq\sqrt{2}$, then $w(S_a)=m(S_a)=1.$
\item[$(b)$]If $\sqrt{2}<a$, then $w(S_a)=\frac{a^2}{2\sqrt{a^2-1}}>m(S_a)=1.$
\end{itemize}
\end{exm}

\begin{prop}\label{periodic}
If $(\omega_n)_{n\geq0}$ is a periodic sequence, then $S$
possesses Bishop's property $(\beta)$.
\end{prop}
\begin{proof}Suppose that $(\omega_n)_{n\geq0}$ is a periodic sequence of period $k$. So, for every $x\in\h$, we have
$\|S^kx\|=\omega_0...\omega_{k-1}\|x\|.$ Hence,
$\frac{1}{\omega_0...\omega_{k-1}}S^k$ is an isometry; so, by
proposition 1.6.7 of \cite{ln}, $S^k$ possesses Bishop's property
$(\beta)$. Therefore, it follows from theorem 3.3.9 of \cite{ln}
that $S$ possesses Bishop's property $(\beta)$.
\end{proof}
Recall that the weight $(\omega_n)_{n\geq0}$ is said to be {\it
almost periodic} if there is a periodic positive sequence
$(p_n)_{n\geq0}$ such that
$\lim\limits_{n\to+\infty}(\omega_n-p_n)=0.$ Note that, if $k$ is
the period of $(p_n)_{n\geq0}$, then
$$r_1(S)=r(S)=(p_0...p_{k-1})^{\frac{1}{k}}.$$ This suggests a
weaker version of question \ref{ques2}.
\begin{ques}\label{ques3}
Suppose that the weight $(\omega_n)_{n\geq0}$ is almost periodic.
Does $S$ possess Bishop's property $(\beta)$?
\end{ques}

\section{Local spectral properties of bilateral weighted shift operators}

In dealing with bilateral weighted shift operators, let
$(\beta_n)_{n\in\Z}$ be the sequence given by

$$ \beta_n=\left\{
\begin{array}{lllll}
\omega_0...\omega_{n-1}&\mbox{if }n>0\\
\\
1&\mbox{if }n=0\\
\\
\frac{1}{\omega_n...\omega_{-1}}&\mbox{if }n<0
\end{array}
\right. $$ and set

$$
\begin{array}{llll}
r^-(S)=\lim\limits_{n\to+\infty}\big{[}\sup\limits_{k>0}
\frac{\beta_{-k}}{\beta_{-n-k}}\big{]}^{\frac{1}{n}},&
r^+(S)=\lim\limits_{n\to+\infty}\big{[}\sup\limits_{k\geq0}
\frac{\beta_{n+k}}{\beta_{k}}\big{]}^{\frac{1}{n}},\\
r_1^-(S)=\lim\limits_{n\to+\infty}\big{[}\inf\limits_{k>0}
\frac{\beta_{-k}}{\beta_{-n-k}}\big{]}^{\frac{1}{n}},&
r_1^+(S)=\lim\limits_{n\to+\infty}\big{[}\inf\limits_{k\geq0}
\frac{\beta_{n+k}}{\beta_{k}}\big{]}^{\frac{1}{n}},\\
r_2^-(S)=\liminf\limits_{n\to+\infty}\big{[}\frac{1}{\beta_{-n}}
\big{]}^{\frac{1}{n}},&r_2^+(S)=
\liminf\limits_{n\to+\infty}\big{[}\beta_n\big{]}^{\frac{1}{n}},\\
r_3^-(S)=\limsup\limits_{n\to+\infty}\big{[}
\frac{1}{\beta_{-n}}\big{]}^{\frac{1}{n}},\mbox{ and }&
r_3^+(S)=\limsup\limits_{n\to+\infty}\big{[}\beta_n\big{]}^{\frac{1}{n}}.
\end{array}
$$ Note that, $$\begin{array}{ll}
r_1(S)=\min(r_1^-(S),r_1^+(S)),&r(S)=\max(r^-(S),r^+(S)),\\
\\
r_1^-(S)\leq r_2^-(S)\leq r_3^-(S)\leq r^-(S),\mbox{ and}&r_1^+(S)
\leq r_2^+(S)\leq r_3^+(S)\leq r^+(S).
\end{array}
$$

\subsection{On the local spectra of bilateral weighted shift operators}
We begin this section by pointing out that there are bilateral
weighted shift operators which do not have the single-valued
extension property. The following result, which can easily be
proved, gives a necessary and sufficient condition for a bilateral
weighted shift operator $S$ to enjoy this property.

\begin{prop}\label{svep}
$S$ has the single-valued extension property if and only if
$r_2^-(S)\leq r_3^+(S)$.
\end{prop}
\begin{proof}If $r_2^-(S)\leq r_3^+(S)$, then according theorem 9 of \cite{shi},
$\sigma_p(S)$ has empty interior. And so, $S$ has the
single-valued extension property.

Conversely, suppose that $r_3^+(S)<r_2^-(S)$. Let
$$O:=\{\lambda\in\C:r_3^+(S)<|\lambda|<r_2^-(S)\}.$$ It is easy to
see that the following analytic $\h-$valued function on $O$,
$$f(\lambda):=\sum\limits_{n\in\Z}\frac{\beta_n}{\lambda^n}e_n,$$
satisfies the equation $$(S-\lambda)f(\lambda)=0,~\lambda\in O.$$
This shows that $S$ does not have the single-valued extension
property.
\end{proof}
\begin{cor}\label{svep*}
$S^*$ has the single-valued extension property if and only if
$r_2^+(S)\leq r_3^-(S)$.
\end{cor}
\begin{proof}As $S^*$ is also a bilateral weighted shift operator with $r_2^-(S^*)=r_2^+(S)$, and $r_3^+(S^*)=r_3^-(S)$, the
result is a direct consequence of proposition \ref{svep}.
\end{proof}

An analogue of proposition \ref{r_2} is given.
\begin{prop}\label{local}
For every non-zero $x\in\h$, we have
$$\{\lambda\in\C:r_3^-(S)<|\lambda|<
r_2^+(S)\}\subset\sigma_{_S}(x).$$
\end{prop}
\begin{proof}If $r_2^+(S)\leq r_3^-(S)$, then there is nothing to prove.
Assume that $r_3^-(S)<r_2^+(S)$, and set
$$O:=\{\lambda\in\C:r_3^-(S)< |\lambda|<r_2^+(S)\}.$$ Now,
consider the following analytic $\h$-valued function $k$ defined
on $O$ by
$$k(\lambda)=\sum\limits_{n\in\Z}\frac{\lambda^n}{\beta_n}e_n.$$
We have $(S-\lambda)^*k(\overline\lambda)=0$ for every $\lambda\in
O$. And so, the rest of the proof goes as in the proof of
proposition \ref{r_2}.
\end{proof}

\begin{cor}\label{concon}If $r_3^-(S)<r_2^+(S)$, then each $\sigma_{_S}(x)$ is connected.
\end{cor}

\begin{proof}The proof is similar to the proof of corollary \ref{con}.
\end{proof}

\begin{rem}A similar proof to
that of (\ref{rsx}) yields
\begin{equation}\label{rsx2}
r_3^+(S)\leq r_{_S}(x)\leq r(S),~\mbox{ for every non-zero }x\in\h.
\end{equation}
If $S$ is invertible, then $S^{-1}$ is also a bilateral weighted
shift with the corresponding weight
$(\frac{1}{\omega_n})_{n\in\Z}$. So, by applying (\ref{rsx2}) to
$S^{-1}$, we get
\begin{equation}\label{22rsx}
\frac{1}{r_2^-(S)}\leq r_{_{S^{-1}}}(x)\leq r\big{(}S^{-1}\big{)},~\mbox{ for every non-zero }x\in\h.
\end{equation}
\end{rem}

Let $\h_0$ denote all finite linear combinations of the vectors
$e_n,~(n\in\Z)$. The following proposition gives a complete
description of the local spectrum of $S$ at a point $x\in\h_0$.

\begin{prop}\label{h0}For every non-zero $x\in\h_0$, we have
$$\sigma_{_S}(x)=\{\lambda\in\C:r_2^-(S)\leq|\lambda|\leq
r_3^+(S)\}.$$
\end{prop}
\begin{proof}
First, let us prove that
\begin{equation}\label{nz}
\sigma_{_S}(e_n)=\{\lambda\in\C:r_2^-(S)\leq|\lambda|\leq
r_3^+(S)\},~\forall n\in\Z.
\end{equation}
To prove (\ref{nz}), we shall limit ourselves to considering the
case $n=0$; the proof for arbitrary $n\in\Z$ is similar.

Let $O_1:=\{\lambda\in\C:|\lambda|> r_3^+(S)\},$ and consider the
following $\h$--valued function, defined on $O_1$ by
$$f_1(\lambda):=\sum\limits_{n=0}^{+\infty}\frac{\beta_n}{\lambda^{n+1}}e_n.$$
As $(S-\lambda)f_1(\lambda)=e_0$ for every $\lambda\in O_1$, we
note that
\begin{equation}\label{f1}
\sigma_{_S}(e_0)\subset\{\lambda\in\C:|\lambda|\leq r_3^+(S)\}.
\end{equation}
We also have
\begin{equation}\label{f2}
\sigma_{_S}(e_0)\subset\{\lambda\in\C:|\lambda|\geq r_2^-(S)\}.
\end{equation}
Indeed, we may assume $r_2^-(S)>0$. The following $\h$--valued
function, defined on $O_2:=\{\lambda\in\C:|\lambda|<r_2^-(S)\}$ by
$$f_2(\lambda):=\sum\limits_{n=1}^{+\infty}\lambda^{n-1}\beta_{-n}e_{-n},$$
satisfies  $(S-\lambda)f_2(\lambda)=e_0$ for every $\lambda\in
O_2$. This proves (\ref{f2}). And so, it follows from (\ref{f1}),
and (\ref{f2}) that
\begin{equation}\label{inclusion1}
\sigma_{_S}(e_0)\subset\{\lambda\in\C:r_2^-(S)\leq|\lambda|\leq
r_3^+(S)\}.
\end{equation}
If $r_2^-(S)>r_3^+(S)$, then
$$\sigma_{_S}(e_0)=\{\lambda\in\C:r_2^-(S)\leq|\lambda|\leq
r_3^+(S)\}=\emptyset.$$ Thus, we may assume now that $r_2^-(S)\leq
r_3^+(S)$. Recall that the condition $r_2^-(S)\leq r_3^+(S)$ means
that $S$ has the single-valued extension property (see proposition
\ref{svep}). As in the beginning of the proof of theorem
\ref{lsp}, one can show that $\sigma_{_S}(e_0)$ has circular
symmetry about the origin. Since $r_{_S}(e_0)=r_3^+(S)$, we have
\begin{equation}\label{ee00}
\{\lambda\in\C:|\lambda|=r_3^+(S)\}\subset\sigma_{_S}(e_0).
\end{equation}
Now, let us show the following

\begin{equation}\label{inclusion2}
\{\lambda\in\C:r_2^-(S)\leq|\lambda|\leq
r_3^+(S)\}\subset\sigma_{_S}(e_0).
\end{equation}
Indeed, from (\ref{ee00}) we can assume that $r_2^-(S)<r_3^+(S)$.
By writing
$\widetilde{e_0}(\lambda)=\sum\limits_{n\in\Z}A_n(\lambda)e_n$,
$\lambda\in\rho_{_S}(e_0)$, we get from the equation
$$(S-\lambda)\widetilde{e_0}(\lambda)=e_0,~\lambda\in\rho_{_S}(e_0),$$
that for every $\lambda\in\rho_{_S}(e_0),$ we have
\begin{equation}\label{ee2}
A_{-1}(\lambda)\omega_{-1}-\lambda A_0(\lambda)=1,
\end{equation}
and
\begin{equation}\label{ee3}
A_n(\lambda)\omega_n-\lambda A_{n+1}(\lambda)=0,~n\in\Z,~n\not=-1.
\end{equation}
And so, from (\ref{ee3}) we see that for every
$\lambda\in\rho_{_S}(e_0),$ we have
\begin{equation}\label{ee4}
\frac{\lambda^n}{\beta_n}A_n(\lambda)=A_0(\lambda),\mbox{ and
}A_{-n}(\lambda)=\omega_{-1}\beta_{-n}\lambda^{n-1}A_{-1}(\lambda),~n\geq1.
\end{equation}
Now, suppose that there is $\lambda_0\in\rho_{_S}(e_0)$ such that
$r_2^-(S)<|\lambda_0|<r_3^+(S)$. So, there are two subsequences
$(n_k)_{k\geq0}$, and $(m_k)_{k\geq0}$ such that for every
$k\geq0$, we have $$|\lambda_0|^{n_k}\leq\beta_{n_k},\mbox{ and
}\frac{1}{\beta_{-m_k}}\leq|\lambda_0|^{m_k}.$$ Therefore, it
follows from (\ref{ee4}) that for every $k\geq0$, we have
$$|A_0(\lambda_0)|\leq |A_{n_k}(\lambda_0)|,\mbox{ and
}|A_{-{m_k}}(\lambda_0)|\geq\omega_{-1}|A_{-1}(\lambda_0)|.$$
Since
$\lim\limits_{n\to+\infty}A_n(\lambda_0)=\lim\limits_{n\to+\infty}A_{-n}(\lambda_0)=0$,
it follows that $A_0(\lambda_0)=A_{-1}(\lambda_0)=0$, which is a
contradiction to (\ref{ee2}). Thus, the desired inclusion
(\ref{inclusion2}) holds. By combining (\ref{inclusion1}), and
(\ref{inclusion2}) we get
\begin{equation}\label{mama}
\sigma_{_S}(e_0)=\{\lambda\in\C:r_2^-(S)\leq|\lambda|\leq
r_3^+(S)\}.
\end{equation}

Finally, let $x=\sum\limits_{n\in\Z}a_ne_n\in\h_0\backslash\{0\}$,
and let $n_0$ be the smallest integer $n$ for which $a_n\not=0$.
So, there is a non-zero polynomial $p$ such that $x=p(S)e_{n_0}$.
If $r_2^-(S)>r_3^+(S)$, then from the fact that
$\sigma_{_S}(p(S)e_{n_0})\subset\sigma_{_S}(e_{n_0})$ and
(\ref{nz}), we see that
$$\sigma_{_S}(x)=\{\lambda\in\C:r_2^-(S)\leq|\lambda|\leq
r_3^+(S)\}=\emptyset.$$ So, assume now that $r_2^-(S)\leq
r_3^+(S)$. If $\sigma_p(S)=\emptyset$, then by theorem
\ref{williams}, we have
$$\sigma_{_S}(x)=\sigma_{_S}(e_{n_0})=\{\lambda\in\C:r_2^-(S)\leq|\lambda|\leq
r_3^+(S)\}.$$ If $\sigma_p(S)\not=\emptyset$, then, in view of
theorem 9(ii) of \cite{shi} and the fact that $r_2^-(S)\leq
r_3^+(S)$, we have $r_2^-(S)=r_3^+(S)$ and
$\sigma_p(S)=\{\lambda\in\C:|\lambda|=r_3^+(S)\}$. Therefore, by
theorem \ref{williams}, we have
$$\sigma_{_S}(e_{n_0})=\sigma_{_S}(x)\cup Z_{_S}(p,e_{n_0}).$$ As
$Z_{_S}(p,e_{n_0})$ is a finite set, and
$\sigma_{_S}(e_{n_0})=\{\lambda\in\C:|\lambda|=r_3^+(S)\}$, we see
that
$$\sigma_{_S}(x)=\sigma_{_S}(e_{n_0})=\{\lambda\in\C:|\lambda|=r_3^+(S)\}.$$
And the desired conclusion holds.
\end{proof}

For every non-zero $x=\sum\limits_{n\in\Z}a_ne_n\in\h$, we set
$$R_\omega^-(x):=\limsup\limits_{n\to+\infty}\big{|}
\frac{a_{-n}}{\beta_{-n}} \big{|}^{\frac{1}{n}},\mbox{ and
}R_\omega^+(x):=\liminf\limits_{n\to+\infty}\big{|}
\frac{\beta_n}{a_n} \big{|}^{\frac{1}{n}}.$$ Note that $0\leq
R_\omega^-(x)\leq r_3^-(S)\mbox{ and }r_2^+(S)\leq
R_\omega^+(x)\leq+\infty$ for every non-zero $x\in\h$.

\begin{thm}\label{lspbil}Assume that $r_2^-(S)\leq r_3^+(S)$. For every non-zero $x\in\h$, the following statements hold.
\begin{itemize}
\item[$(a)$]If $R_\omega^-(x)<r_2^-(S)\mbox{ and }r_3^+(S)< R_\omega^+(x)$, then
$$\sigma_{_S}(x)=\{\lambda\in\C:r_2^-(S)\leq|\lambda|\leq
r_3^+(S)\}.$$
\item[$(b)$]Otherwise,
$$\{\lambda\in\C:\max(R_\omega^-(x),r_2^-(S))<|\lambda|<\min(R_\omega^+(x),r_3^+(S))\}\subset\sigma_{_S}(x).$$
\end{itemize}
\end{thm}
\begin{proof}Recall again that the assumption $r_2^-(S)\leq r_3^+(S)$ means that $S$ has the single-valued extension property (see
proposition \ref{svep}). Let $x=\sum\limits_{n\in\Z}a_ne_n$ be a
non-zero element of $\h$.

$(a)$ Assume that $R_\omega^-(x)<r_2^-(S),\mbox{ and }r_3^+(S)<
R_\omega^+(x)$. Keep in mind that
$\sigma_{_S}(e_0)=\{\lambda\in\C:r_2^-(S)\leq|\lambda|\leq
r_3^+(S)\},$ and $$ \widetilde{e_0}(\lambda)=\left\{
\begin{array}{lll}
-\sum\limits_{n=0}^{+\infty}\frac{\beta_n}{\lambda^{n+1}}e_n&\mbox{if
}|\lambda|>r_3^+(S)\\
\\
\sum\limits_{n=1}^{+\infty}\lambda^{n-1}\beta_{-n}e_{-n}&\mbox{if
}|\lambda|<r_2^-(S)
\end{array}
\right. $$ Now, consider the following analytic $\h$--valued
function, defined on the open annulus
$\{\lambda\in\C:R_\omega^-(x)<|\lambda|<R_\omega^+(x)\}$ by
$$f(\lambda):=\sum\limits_{n\in\Z}\frac{a_n}{\beta_n}\lambda^n.$$
Let $R^+$ and $R^-$ be two reals such that
$$r_3^+(S)<R^+<R_\omega^+(x),\mbox{ and
}R_\omega^-(x)<R^-<r_2^-(S).$$ Let $\Gamma^+$ (resp. $\Gamma^-$)
be the circle centered at the origin and of radius $R^+$ (resp.
$R^-$). Suppose that $\Gamma^+$ (resp. $\Gamma^-$) is oriented
anti-clockwise (resp. clockwise) direction. We have
\begin{eqnarray*}
f(S,e_0)&=&\frac{-1}{2\pi
i}\oint_{\lambda\in\Gamma^+\cup\Gamma^-}f(\lambda)\widetilde{e_0}(\lambda)d\lambda\\
&=&\sum\limits_{n=0}^{+\infty}\bigg[\frac{1}{2\pi
i}\oint_{\lambda\in\Gamma^+}\frac{f(\lambda)}{\lambda^{n+1}}
d\lambda \bigg]\beta_ne_n\\
&&+\sum\limits_{n=1}^{+\infty}\bigg[\frac{-1}{2\pi
i}\oint_{\lambda\in\Gamma^-}f(\lambda)\lambda^{n-1} d\lambda
\bigg]\beta_{-n}e_{-n}\\ &=&x.
\end{eqnarray*}
And so, as in the last part of the proof of proposition \ref{h0},
we get from theorem \ref{williams} that
$$\sigma_{_S}(x)=\sigma_{_S}(e_0)=\{\lambda\in\C:r_2^-(S)\leq|\lambda|\leq
r_3^+(S)\}.$$

$(b)$ Suppose that
$\max(R_\omega^-(x),r_2^-(S))<\min(R_\omega^+(x),r_3^+(S))$,
otherwise there is nothing to prove. By writing
$\widetilde{x}(\lambda)=\sum\limits_{n\in\Z}A_n(\lambda)e_n$,
$\lambda\in\rho_{_S}(x)$, we get from the equation
$$(S-\lambda)\widetilde{x}(\lambda)=x,~\lambda\in\rho_{_S}(x),$$
that for every $\lambda\in\rho_{_S}(x)$, we have
\begin{equation}\label{bbbb1}
A_n(\lambda)\omega_n-\lambda A_{n+1}(\lambda)=a_{n+1},~\forall
n\in\Z.
\end{equation}
In particular, for every $\lambda\in\rho_{_S}(x)$, we have
\begin{equation}\label{bbbb2}
A_{-1}(\lambda)\omega_{-1}-\lambda A_0(\lambda)=a_0.
\end{equation}
And so, by a simple calculation, we get from (\ref{bbbb1}) that,
for every $\lambda\in\rho_{_S}(x)$ and $n\geq1$,
\begin{equation}\label{pside}
\frac{\lambda^{n}}{\beta_n}A_n(\lambda)=
A_0(\lambda)-\frac{a_1}{\beta_1}-\frac{a_2}{\beta_2}\lambda-...-\frac{a_n}{\beta_n}\lambda^{n-1},
\end{equation}
and
\begin{equation}\label{nsidep}
A_{-n}(\lambda)=\frac{\lambda^{n-1}}{\omega_{-2}...\omega_{-n}}A_{-1}(\lambda)+\frac{\lambda^{n-2}a_{-1}}{\omega_{-2}...\omega_{-n}}
+...+\frac{a_{-n+1}}{\omega_{-n}}.
\end{equation}
Note that for every $\lambda\in\rho_{_S}(x)$, $\lambda\not=0$,
(\ref{nsidep}) can be reformulated as follows.
\begin{equation}\label{nside}
A_{-n}(\lambda)=\beta_{-n}\lambda^{n-1}\bigg[\omega_{-1}A_{-1}(\lambda)+\frac{a_{-1}}{\beta_{-1}\lambda}+
\frac{a_{-2}}{\beta_{-2}\lambda^2}+...+\frac{a_{-n+1}}{\beta_{-n+1}\lambda^{n-1}}\bigg].
\end{equation}
Now, suppose that
$$O:=\{\lambda\in\C:\max(R_\omega^-(x),r_2^-(S))<|\lambda|<\min(R_\omega^+(x),r_3^+(S))\}\cap\rho_{_S}(x)\not=\emptyset.$$
Fix $\lambda\in O$. There are two subsequences $(n_k)_{k\geq0}$,
and $(m_k)_{k\geq0}$ such that for every $k\geq0$, we have
$$|\lambda|^{n_k}\leq\beta_{n_k},\mbox{ and
}\frac{1}{\beta_{-m_k}}\leq|\lambda|^{m_k}.$$ And so, from
(\ref{pside}) and (\ref{nside}), we get, respectively,
$$|A_{n_k}(\lambda)|\geq|A_0(\lambda)-\frac{a_1}{\beta_1}-\frac{a_2}{\beta_2}\lambda-...-\frac{a_{n_k}}{\beta_{n_k}}\lambda^{n_k-1}|,$$
and
$$|A_{-{m_k}}(\lambda)|\geq|\omega_{-1}A_{-1}(\lambda)+\frac{a_{-1}}{\beta_{-1}\lambda}+
\frac{a_{-2}}{\beta_{-2}\lambda^2}+...+\frac{a_{-{m_k}+1}}{\beta_{-{m_k}+1}\lambda^{{m_k}-1}}|.$$
As
$\lim\limits_{n\to+\infty}A_n(\lambda)=\lim\limits_{n\to+\infty}A_{-n}(\lambda)=0$,
and both series
$\sum\limits_{n\geq1}\frac{a_n}{\beta_n}\lambda^{n-1}$ and
$\sum\limits_{n\leq-1}\frac{a_n}{\beta_n}\lambda^n$ converge, it
follows that
\begin{equation}\label{ppside}
A_0(\lambda)=\sum\limits_{n=1}^{+\infty}\frac{a_n}{\beta_n}\lambda^{n-1},
\end{equation}
and
\begin{equation}\label{nnside}
\omega_{-1}A_{-1}(\lambda)=-\sum\limits_{n\leq-1}\frac{a_n}{\beta_n}\lambda^n.
\end{equation}
So, from (\ref{bbbb2}), (\ref{ppside}), and (\ref{nnside}), we see
that $\sum\limits_{n\in\Z}\frac{a_n}{\beta_n}\lambda^n=0$ for
every $\lambda\in O$; therefore, $a_n=0$ for every $n\in\Z$. This
is a contradiction to the fact that $x\not=0$. And, the proof is
complete.
\end{proof}

Note that the importance of theorem \ref{lspbil} occurs when $r_2^-(S)>0$. In this case, we computed
the local spectrum of $S$ at most points in $\h$. The next result deals with the case $r_2^-(S)=0$.
For every $k\in\Z$, we write
$$\h_k^+=\bigvee\{e_n:n\geq k\}\mbox{ and }\h_k^-=\bigvee\{e_n:n<k\},$$
where "$\bigvee$" denotes the closed linear span.

\begin{thm}\label{particular}Assume that $r_2^-(S)=0$. Let $k\in\Z$, then for every non-zero $x\in\h_k^+$, the following statements hold.
\begin{itemize}
\item[$(a)$]If $r_3^+(S)< R_\omega^+(x)$, then
$\sigma_{_S}(x)=\{\lambda\in\C:|\lambda|\leq r_3^+(S)\}.$
\item[$(b)$]If $r_3^+(S)\geq R_\omega^+(x)$, then
$\{\lambda\in\C:|\lambda|\leq R_\omega^+(x)\}\subset\sigma_{_S}(x).$
\end{itemize}
\end{thm}
\begin{proof}
Since $\h_k^+$ and $\h_k^-$ are invariant subspaces of $S$ and $S^*$ respectively, we note that
$S^+:=S_{|\h_k^+}$ and $S^-:=S^*_{|\h_k^-}$ are
unilateral weighted shift operators with corresponding weight sequences
$\omega^+:=(\omega_n)_{n\geq k}$ and $\omega^-:=(\omega_n)_{n<k}$. We have
$$r_i(S^\pm)=r_i^\pm(S),\mbox{ and }r(S^\pm)=r^\pm(S),~~(i=1,2,3).$$
Moreover, for every $x\in\h_k^+$, we have
$$R_\omega^+(x)=R_{\omega^+}(x).$$
We note also that, since $r_2(S^-)=r_2^-(S)=0$,
${S^-}^*$ has the single-valued extension property. As $S$ has the following matrix
representation on $\h=\h_k^+\oplus\h_k^-$:
$$S=\left[
\begin{array}{cc}
S^+&\ast\\
0&{S^-}^*
\end{array}
\right],$$
it follows from lemma \ref{kim33} that $\sigma_{_S}(x)=\sigma_{_{S^+}}(x)$
for every $x\in\h_k^+$. In view of theorem \ref{lsp}, the proof is complete.
\end{proof}

\subsection{Which bilateral weighted shift operators satisfy Dunford's condition $(C)$?}
The following result gives a necessary condition for the bilateral
weighted shift operator $S$ to satisfy Dunford's condition $(C)$.

\begin{thm}\label{dccbil}Assume that $S$ satisfies Dunford's condition $(C)$.
\begin{itemize}
\item[$(a)$] If $S$ is not invertible, then
\begin{equation}\label{dcc2}
r_1(S)=r_2^-(S)=0\leq r_3^+(S)=r(S).
\end{equation}
\item[$(b)$] If $S$ is invertible, then
\begin{equation}\label{dcc1}
r_1(S)=\frac{1}{r(S^{-1})}=r_2^-(S)\leq r_3^+(S)=r(S).
\end{equation}
\end{itemize}
\end{thm}
\begin{proof}
Note that, since $S$ satisfies Dunford's condition $(C)$, $S$ has
the single-valued extension property. By proposition \ref{svep},
we have $r_2^-(S)\leq r_3^+(S)$. Now, let
$F:=\{\lambda\in\C:r_2^-(S)\leq|\lambda|\leq r_3^+(S)\}$; it
follows from proposition \ref{h0} that $\h_{_S}(F)$ contains a
dense subspace. As $\h_{_S}(F)$ is a closed subspace, we have
$\h_{_S}(F)=\h$; so, $\sigma_{_S}(x)\subset F$ for every $x\in\h$.
Since, $\sigma(S)=\bigcup\limits_{x\in\h}\sigma_{_S}(x)\subset F$
(see proposition 1.3.2 of \cite{ln}), the desired conclusion holds
from theorem 5 of \cite{shi}.
\end{proof}

\begin{rem}\label{view}
An interesting special case occurs when the sequences
$([\beta_n]^{\frac{1}{n}})_{n\geq1}$ and $\big(\big[\frac{1}{\beta_{-n}}\big]\big)_{n\geq1}$
converge respectively to $a$ and $b$.
By combining proposition \ref{local}, and
 theorem \ref{dccbil}, we see that

1. If $S$ is not invertible, then the following statement are equivalent.
\begin{itemize}
\item[$(a)$] $S$ has fat local spectra, that is $\sigma_{_S}(x)=\sigma(S)$
for every non-zero $x\in\h$.
\item[$(b)$] $S$ satisfies Dunford's condition $(C)$.
\item[$(c)$] $b=0\leq a=r(S)$.
\end{itemize}

2. If $S$ is invertible and its spectrum is not a circle, then following
statement are equivalent.
\begin{itemize}
\item[$(a)$] $S$ has fat local spectra.
\item[$(b)$] $S$ satisfies Dunford's condition $(C)$.
\item[$(c)$] $b=\frac{1}{r(S^{-1})}<a=r(S)$.
\end{itemize}
\end{rem}

As the unilateral case, we note that in general setting if
$S$ is not invertible (resp. invertible) and satisfies (\ref{dcc2}) (resp.
(\ref{dcc1})), then $\sigma_{_S}(x)=\sigma(S)$
for all non-zero $x$ in a dense subset of $\h$ (see proposition \ref{h0},
theorem \ref{lspbil} and theorem \ref{particular}). Thus, we conjecture that

1. If $S$ is not invertible, then the following statements are equivalent.
\begin{itemize}
\item[$(a)$] $S$ has fat local spectra.
\item[$(b)$] $S$ satisfies Dunford's condition $(C)$.
\item[$(c)$] $r_1(S)=r_2^-(S)=0\leq r_3^+(S)=r(S)$.
\end{itemize}

2. If $S$ is invertible and its spectrum is not a circle, then the following statements are equivalent.
\begin{itemize}
\item[$(a)$] $S$ has fat local spectra.
\item[$(b)$] $S$ satisfies Dunford's condition $(C)$.
\item[$(c)$] $r_1(S)=r_2^-(S)$, and $r_3^+(S)=r(S)$.
\end{itemize}

Note that if the spectrum of a bilateral weighted shift operator
$S$ is a circle, then $S$ may have fat local spectra
(see example \ref{atzmon}) as well $S$ may satisfy
Dunford's condition $(C)$ and without fat local spectra. This is the case
for example of any non quasi-nilpotent decomposable bilateral weighted
shift operator (see corollary \ref{deco}). But, we do not
know in general, if a bilateral weighted shift operator $S$ satisfies
always Dunford's condition $(C)$ whenever its spectrum is a circle.

\subsection{Which bilateral weighted shift operators possess Bishop's property $(\beta)$?}

For an operator $T\in\lh$, we shall denote
$$q(T):=\min\{|\lambda|:\lambda\in\sigma(T)\}.$$
Note that
$$q(T)=\left\{
\begin{array}{lll}
0&\mbox{ if }T\mbox{ is not invertible}\\
\\
\frac{1}{r(T^{-1})}&\mbox{ if }T\mbox{ is invertible}.
\end{array}
\right.$$

Now, suppose that $S$ possesses Bishop's property $(\beta)$. It is shown in
theorem 3.10 of \cite{bourhim} that

$(a)$ If $S$ is invertible then $S$ satisfies both (\ref{qs1}) and (\ref{qs2}).

$(b)$ If $S$ is not invertible then $S$ satisfies (\ref{qs2}) and
$q(S)=r_2^-(S)=0$.

\noindent On the other hand, it is shown in theorem 2.7 of \cite{mmn} that
$$q(S)=r_1(S)=r_2^-(S)=r_3^-(S),$$ and $$r_2^+(S)=r_3^+(S)=r^+(S)=r(S).$$ Therefore, the authors
deduce from proposition \ref{local} that either $S$ has fat local spectra or $r_1(S)=r(S)$.
Using some ideas from both proofs given in \cite{bourhim} and \cite{mmn}, we refine these results as follows.

\begin{thm}\label{bila}If $S$ possesses Bishop's property $(\beta)$, then
\begin{equation}\label{qs1}
q(S)=r_1^-(S)=r_2^-(S)=r_3^-(S)=r^-(S),
\end{equation}
and
\begin{equation}\label{qs2}
r_1^+(S)=r_2^+(S)=r_3^+(S)=r^+(S)=r(S).
\end{equation}
\end{thm}

In order to prove this theorem, we shall use the following result from \cite{mmn}.

\begin{prop}\label{localmmn}
Assume that $T\in\lh$ is an injective operator. If $T$ possesses Bishop's
property $(\beta)$ then for every $x\in T^\infty\h$, the sequence $(\|T^{-n}x\|^{-\frac{1}{n}})_{n\geq1}$ is
convergent and
$\lim\limits_{n\to+\infty}\|T^{-n}x\|^{-\frac{1}{n}}=\min\{|\lambda|:\lambda\in\sigma_{_T}(x)\}.$
\end{prop}
\noindent{\it Proof of theorem }\ref{bila}. Let $\h_0^+:=\bigvee\{e_k:k\geq0\}$ and let
$S^+:=S_{|\h_0^+}$. It is clear that
$S^{+}$ is a unilateral weighted shift operator with Bishop's
property $(\beta)$. Hence, according theorem \ref{beta1}, we have
$r_1(S^{+})=r(S^{+})$. As
$$r_i(S^{+})=r_i^+(S),~r(S^{+})=r^+(S),~(i=1,2,3),$$
and
$$r_1^+(S)\leq r_2^+(S)\leq r_3^+(S)\leq r^+(S),$$
we have
\begin{equation}\label{bbil1}
r_1^+(S)=r_2^+(S)=r_3^+(S)=r^+(S).
\end{equation}
Note that, since $S$ satisfies Dunford's condition $(C)$, we have
$r_3^+(S)=r(S)$ (see theorem \ref{dccbil}).
Therefore, from (\ref{bbil1}), we see that
$$r_1^+(S)=r_2^+(S)=r_3^+(S)=r^+(S)=r(S).$$
Since $r(S)=\max(r^-(S),r^+(S))$, we have
\begin{equation}\label{bbil11}
r^-(S)\leq r_1^+(S)=r_2^+(S)=r_3^+(S)=r^+(S)=r(S).
\end{equation}

If $S$ is invertible, then $S^{-1}$ possesses also Bishop's property
$(\beta)$ (see theorem 3.3.9 of \cite{ln}). As,
$$r^\pm(S^{-1})=\frac{1}{r_1^\mp(S)},~r_1^\pm(S^{-1})=\frac{1}{r^\mp(S)}\mbox{ and }r_{2,3}^\pm(S^{-1})=\frac{1}{r_{3,2}^\mp(S)},
$$
the desired conclusion holds
by applying (\ref{bbil11}) to $S^{-1}$.

Now, assume that $S$ is not invertible. Since $S$ satisfies Dunford's condition $(C)$, it follows from (\ref{dcc2}) that
$r_1(S)=r_2^-(S)=0$. Hence,
\begin{equation}\label{ration1}
q(S)=r_1(S)=r_1^-(S)=r_2^-(S)=0.
\end{equation}
In view of (\ref{bbil11}) and (\ref{ration1}), we get from theorem 7 of \cite{shi} that
$$\sigma_{ap}(S)=\{\lambda\in\C:|\lambda|\leq r^-(S)\}\cup\{\lambda\in\C:|\lambda|=r(S)\}.$$
Since $S$ is a rationally cyclic injective operator, it follows from proposition \ref{rationally} that
\begin{eqnarray*}
\Re(S^*)&=&\overline{\sigma(S)\backslash\sigma_{ap}(S)}\\
&=&\{\lambda\in\C:r^-(S)<|\lambda|<r(S)\}.
\end{eqnarray*}
On the other hand, we always have
$$\Re(S^*)=\{\lambda\in\C:r_3^-(S)<|\lambda|<r_2^+(S)\}.$$
Therefore,
\begin{equation}\label{ration2}
r_3^-(S)=r^-(S).
\end{equation}
Since the sequence $(\|S^{-n}e_0\|^{-\frac{1}{n}})_{n\geq1}$ is convergent (see proposition \ref{localmmn}), we have
\begin{equation}\label{ration3}
r_2^-(S)=r_3^-(S).
\end{equation}
Thus, the desired result follows from (\ref{ration1}), (\ref{ration2}) and (\ref{ration3}). $\hfill\square$

Conversely, we have the following result.

\begin{thm}Assume that $S$ satisfies $(\ref{qs1})$ and $(\ref{qs2})$. The following statements hold.
\begin{itemize}
\item[$(a)$]If $S$ is not invertible, then
either $S$ possesses Bishop's property $(\beta)$ or
$\sigma_\beta(S)=\{\lambda\in\C:|\lambda|=r(S)\}.$
\item[$(b)$]If $S$ is invertible, then either $S$ possesses Bishop's property $(\beta)$ or $\sigma_\beta(S)$ equals either one of the circles
$\{\lambda\in\C:|\lambda|=\frac{1}{r(S^{-1})}\}$ and $\{\lambda\in\C:|\lambda|=r(S)\}$
or the union of both of them.
\end{itemize}
\end{thm}
\begin{proof}
$(a)$ Assume that $S$ is not invertible. If we set $\h_0^+=\bigvee\{e_k:k\geq0\}$ and
$\h_0^-=\bigvee\{e_k:k<0\}$, then $S$ has the following matrix
representation on $\h=\h_0^+\oplus\h_0^-$:
$$S=\left[
\begin{array}{cc}
S^+&*\\
0&{S^-}^*
\end{array}
\right],$$
where $S^+=S_{|\h_0^+}$ and $S^-=S^*_{|\h_0^-}$ are
unilateral weighted shift operators with corresponding weight sequences
$\omega^+:=(\omega_n)_{n\geq 0}$ and $\omega^-:=(\omega_n)_{n<0}$.
Since $r(S^-)=r^-(S)=0$ (see (\ref{qs1})), $S^-$ is quasi-nilpotent. In particular, ${S^-}^*$
possesses Bishop's property $(\beta)$. As
$r_1(S^+)=r_1^+(S)$, $r(S^+)=r^+(S)$ and $r_1^+(S)=r^+(S)=r(S)$ (see (\ref{qs2})), we have
$r_1(S^+)=r(S^+)=r(S)$. By theorem \ref{beta1}, we see that either
$S^+$ possesses Bishop's property $(\beta)$ or $\sigma_\beta(S^+)=\{\lambda\in\C:|\lambda|=r^+(S)\}$.
Now, the desired conclusion follows from lemma \ref{kim2}.

$(b)$ Assume that $S$ is invertible and satisfies (\ref{qs1}) and (\ref{qs2}). It follows from theorem 7 of \cite{shi}
that $$\sigma_{ap}(S)=\{\lambda\in\C:|\lambda|=\frac{1}{r(S^{-1})}\}\cup\{\lambda\in\C:|\lambda|=r(S)\}.$$
By proposition \ref{kim1}, and lemma \ref{kim5}, the result is
established.
\end{proof}

\begin{rem}As in the unilateral case, we note that if $r_1(S)=r(S)$, then $S$ possesses
Bishop's property $(\beta)$ or
$\sigma_\beta(S)=\{\lambda\in\C:|\lambda|=r(S)\}$.
But, if $S$ possesses Bishop's property $(\beta)$, then $r_1(S)$ and
$r(S)$ need not be equal as the following example shows. Indeed,
let $S$ be the bilateral weighted shift with the corresponding
weight $(\omega_n)_{n\in\Z}$ given by $$ \omega_n=\left\{
\begin{array}{lll}
2&\mbox{if }n\geq0\\
\\
1&\mbox{if }n<0
\end{array}
\right. $$ It is clear that $S$ is a hyponormal operator; in
particular, $S$ possesses Bishop's property $(\beta)$. On the
other hand, we have $ r_1(S)=1<r(S)= 2.$
\end{rem}

\begin{rem}
In view of remark \ref{view}, and theorem \ref{bila}, one can
construct examples of bilateral weighted shift operators without
Bishop's property $(\beta)$ which satisfy Dunford's condition $(C)$
(see also example \ref{atzmon}).
\end{rem}

We give here a simple proof of the well known following result.

\begin{cor}\label{deco}
If $S$ is decomposable then either $S$ is quasi-nilpotent or it is
invertible and its spectrum is a circle. In particular,
$r_1(S)=r(S)$.
\end{cor}
\begin{proof}If $S$ is decomposable, then both $S$ and $S^*$ possess
Bishop's property $(\beta)$. Since $S^*$ is also a bilateral weighted
shift and satifies $$r_i^\pm(S^*)=r_i^\mp(S)\mbox{ and }r^\pm(S^*)=r^\mp(S),~~(i=1,2,3),$$
the desired result holds by applying theorem \ref{bila} to $S$ and $S^*$.
\end{proof}

The next two propositions give a sufficient condition on the weight
$(\omega_n)_{n\in\Z}$ so that the corresponding bilateral weighted shift
$S$ is decomposable.

\begin{prop}
If $m(S)=w(S)$, then $S$ is decomposable.
\end{prop}
\begin{proof}Without loss of generality, assume that $m(S)=w(S)=1$.
As in the proof of proposition \ref{w}, for every
$\lambda\not\in\sigma(S)=\{\mu\in\C:|\mu|=1\}$
and for every $x\in\h$, we have $$|1-|\lambda||\|x\|\leq
\|(S-\lambda)x\|.$$ Therefore, it follows from theorem 1.7.2 of
\cite{ln} that $S$ is decomposable.
\end{proof}

\begin{prop}
If $(\omega_n)_{n\in\Z}$ is a periodic sequence, then $S$ is
decomposable.
\end{prop}
\begin{proof}Suppose that $(\omega_n)_{n\in\Z}$ is a periodic sequence of
period $k$. So, as in the proof of proposition
\ref{periodic}, we have $\frac{1}{\omega_0...\omega_{k-1}}S^k$ is
an isometry which is clearly invertible. By, proposition 1.6.7 of
\cite{ln}, $S^k$ is decomposable. And so, it follows from theorem
3.3.9 of \cite{ln} that $S$ is decomposable.
\end{proof}
Recall that the weight $(\omega_n)_{n\in\Z}$ is said to be {\it
almost periodic} if there is a periodic positive sequence
$(p_n)_{n\in\Z}$ such that
$\lim\limits_{|n|\to+\infty}(\omega_n-p_n)=0.$ Note that, if $k$
is the period of $(p_n)_{n\in\Z}$, then
$$r_1(S)=r(S)=(p_0...p_{k-1})^{\frac{1}{k}}.$$
So, one may ask if $S$ is decomposable when the weight $(\omega_n)_{n\in\Z}$ is almost periodic.
It turn out that this is not true as the following example shows.

\begin{exm}\label{atzmon}Consider the following real-valued
continuous function on $\R$
$$\varphi(x):=\exp\bigg( \frac{|x|\sin\big[\log\log\log(|x|+e^3)\big]}
{\log(|x|+e)}\bigg),$$
and set
$$\omega_n=\exp\big(\varphi(n+1)-\varphi(n)\big),~~n\in\Z.$$
Since $\lim\limits_{|x|\to+\infty}\varphi^{'}(x)=0$, we have
$$\lim\limits_{|n|\to+\infty}\omega_n=1.$$
Hence, $(\omega_n)_{n\in\Z}$ is in particular almost periodic and
$$r_1(S)=r(S)=1.$$
On the other hand, it shown in \cite{atzmon} that $(\omega_n)_{n\in\Z}$
satisfies all the hypotheses of theorem 2 of \cite{atzmon}. Therefore,
it follows from the proof of theorem 2 of \cite{atzmon} that
$$\sigma_{_S}(x)=\sigma_{_{S^*}}(x)=\sigma(S)=\{\lambda\in\C:|\lambda|=1\}
\mbox{ for every non-zero }x\in\h.$$
This shows that both $S$ and $S^*$ have hat local spectra, and then
satisfy Dunford's condition $(C)$,
but neither $S$ nor $S^*$ possesses Bishop's property $(\beta)$.
\end{exm}

\section{Examples and comments}

The {\it quasi-nilpotent part} of an operator $T\in\lh$ is the set
$$\h_0(T)=\{x\in\h:\lim\limits_{n\to+\infty}\|T^nx\|^{\frac{1}{n}}=0\}.$$
It is a linear subspace of $\h$, generally not closed. The
operator $T$ is said to have {\it property } $(Q)$ if
$\h_0(T-\lambda)$ is closed for every $\lambda\in\C$. Note that
Dunford's condition $(C)$ implies property $(Q)$ and in its turn
property $(Q)$ implies the single-valued extension property. In
\cite{aiena}, the authors have shown, by a convolution operator of
group algebras and a direct sum of unilateral weighted shits, that
the above implications are not reversed in general. The purpose of
this section is to provide elementrary results in order to produce
simpler counter-examples showing that property $(Q)$ is strictly intermediate to
Dunford's condition $(C)$ and the single-valued extension
property.

Throughout the remainder of this paper, we suppose that $S$ is a
unilateral weighted shift operator on $\h$ and the weights
$\omega_n,~(n\geq0)$, are non-negative. Note that $S$ is injective
if and only if none of the weights is zero.

It is shown in proposition 17 of \cite{amn} that if $S$ is an injective
unilateral weighted shift operator, then $\h_0(S)+\ran(S)$ is dense in $\h$
if and only if $r_3(S)=0$. However, in view of (\ref{rsx}) and the fact that
$r_{_S}(e_n)=r_3(S)$ for all $n\geq0$, one can see that either $\h_0(S)=\{0\}$
or $\h_0(S)$ is dense in $\h$. More precisely, we have
\begin{prop}\label{jawad}If $S$ is an injective unilateral weighted shift operator, then
$\h_0(S)$ is dense in $\h$ if and only if $r_3(S)=0$.
\end{prop}

The following two propositions enable us to produce examples of
operators which have the single-valued extension property and
without property $(Q)$.

\begin{prop}\label{sstar}The following assertions are equivalent.
\begin{itemize}
\item[$(a)$]$S^*$ is quasi-nilpotent.
\item[$(b)$]$S^*$ is decomposable.
\item[$(c)$]$S^*$ possesses Bishop's property $(\beta)$.
\item[$(d)$]$S^*$ satisfies Dunford's condition $(C)$.
\item[$(e)$]$S^*$ has property $(Q)$.
\item[$(f)$]$\h_0(S^*)$ is closed.
\end{itemize}
\end{prop}
\begin{proof}
It suffices to show that the implication $(f)\Rightarrow(a)$
holds. For every $n\geq0$, we have ${S^*}^{^{n+1}}e_n=0$. This
shows that $e_n\in\h_0(S^*)$ for every $n\geq0$. Since $\h_0(S^*)$
is closed, we have $\h_0(S^*)=\h.$ And so, it follows from theorem
1.5 of \cite{praha} that $S^*$ is quasi-nilpotent.
\end{proof}

A similar proof of proposition \ref{sstar} yields the following
result.

\begin{prop}\label{star}If infinitely many weights $\omega_n$ are zero, then the following statments are equivalent.
\begin{itemize}
\item[$(a)$]$S$ is quasi-nilpotent.
\item[$(b)$]$S$ is decomposable.
\item[$(c)$]$S$ possesses Bishop's property $(\beta)$.
\item[$(d)$]$S$ satisfies Dunford's condition $(C)$.
\item[$(e)$]$S$ has property $(Q)$.
\item[$(f)$]$\h_0(S)$ is closed.
\end{itemize}
\end{prop}

\begin{exm}Note that if $S$ is injective, then $S^*$ has the single-valued
extension property if and only if $r_2(S)=0$. Thus, in view of
proposition \ref{sstar}, we see that if $0=r_2(S)<r(S)$, then
$S^*$ has the single-valued extension property but not property
$(Q)$.

It remains to produce a positive bounded weight sequence
$(\omega_n)_{n\geq0}$ such that the corresponding weighted shift
operator $S$ satisfies $0=r_2(S)<r(S)$. Let $(k_i)_{i\geq1}$ be
the sequence given by $$k_1=1\mbox{ and
}k_{i+1}=(i+1)k_i+1,~\forall i\geq1.$$ Let $(\omega_n)_{n\geq0}$
be the weight given by $$\omega_n=\left\{
\begin{array}{lll}
\frac{1}{2^{ik_i}}&\mbox{if }n=k_i-1\mbox{ for some }i\geq1\\
\\
2&\mbox{otherwise}
\end{array}
\right. $$ It is easy to check that $\beta_n\leq1$ for every
$n\geq0$, and that $$\beta_{k_i-1}=1,\mbox{ and
}\beta_{k_i}=\frac{1}{2^{ik_i}}\mbox{ for every }i\geq1.$$ This
shows that $0=r_2(S)<1=r_3(S)\leq r(S)$. Therefore, $S^*$ has the
single-valued extension property and without property $(Q)$.
\end{exm}

\begin{exm}Let $(\omega_n)_{n\geq0}$ be the weight given by
$$\omega_n=\left\{
\begin{array}{lll}
0&\mbox{if }n\mbox{ is a square of an integer}\\
\\
1&\mbox{otherwise}
\end{array}
\right. $$ It is easy to see that $\|S^n\|=1$ for every $n\geq 1$.
This shows that $S$ is not quasi-nilpotent. And so, $S$ has not
property $(Q)$ (see Proposition \ref{star}). On the other hand,
$S$ has the single-valued extension property since
$\sigma_p(S)=\{0\}$.
\end{exm}

If infinitely many weights $\omega_n$ are zero, then it is seen in
proposition \ref{star} that $S$ enjoy property $(Q)$ if and only
if it is quasi-nilpotent. When only finitely many weights
$\omega_n$ are zero, the corresponding weighted shift $S$ is a
direct sum of finite-dimensional nilpotent operator $S_1$ and an
injective unilateral weighted shift operator $S_2$. Therefore, $S$
has property $(Q)$ if and only if $S_2$ has property $(Q)$. Thus
in studying which unilateral weighted shift operator enjoy
property $(Q)$ we may assume that none of the weights $\omega_n$
is zero.

\begin{prop}\label{propertyQ}If $S$ is injective, then the following statments are equivalent.
\begin{itemize}
\item[$(a)$]$S$ has property $(Q)$.
\item[$(b)$]Either $S$ is quasi-nilpotent or $r_3(S)>0$.
\end{itemize}
\end{prop}
\begin{proof}Since $S$ has the single-valued extension property, we have
$\h_0(S-\lambda)=\h_{_S}(\{\lambda\})$ for every
$\lambda\in\C$.

Assume that $r_3(S)>0$. For every non-zero $x\in\h$, we have
$\sigma_{_S}(x)$ is a connected set containing $0$ (see corollary
\ref{con}) and $r_3(S)\leq r_{_S}(x)\leq r(S)$ (see (\ref{rsx})).
This implies that $\h_0(S-\lambda)=\h_{_S}(\{\lambda\})=\{0\}$ for
every $\lambda\in\C$. Therefore, $S$ has property $(Q)$ and the
implication $(b)\Rightarrow(a)$ holds.

Now, suppose that $S$ has property $(Q)$ and $r_3(S)=0$. From
proposition \ref{jawad} we see that  $\h_0(S)$ is a dense subspace of
$\h$. As $\h_0(S)$ is closed, we have $\h_0(S)=\h$. By theorem 1.5
of \cite{praha}, $S$ is quasi-nilpotent. This finishes the proof.
\end{proof}
To separate Dunford's condition $(C)$ and property $(Q)$ we need
to produce a positive bounded weight sequence
$(\omega_n)_{n\geq0}$ such that the corresponding weighted shift
operator $S$ satisfies $0<r_3(S)<r(S)$ (see theorem \ref{lah} and
proposition \ref{propertyQ}).

\begin{exm}
Let $(C_k)_{k\geq 0}$ be a sequence of successive disjoint
segments covering the set $\N$ of non-negative integers such that
each segment $C_k$ contains $k^2$ elements. Let $k\in\N$, for
$n\in C_k$ we set $$ \omega_n=\left\{
\begin{array}{ll}
2&\mbox{if }n\mbox{ is one of the first }k^{th}\mbox{ terms of
}C_k\\
\\
1&\mbox{otherwise }
\end{array}
\right. $$ It is easy to see that $r_1(S)=r_3(S)=1<r(S)=2.$
Therefore, $S$ has property $(Q)$ but not Dunford's condition
$(C)$.

The original idea of this construction is due to W. C. Ridge
\cite{rid}.
\end{exm}

Finally, we would like to point out that

$(a)$ proposition \ref{jawad} remain valid
for an injective bilateral weighted shift operator $S$ by replacing $r_3(S)$
with $r_3^+(S)$.

$(b)$ if $S$ is a non injective bilateral weighted shift operator then $S$ is a
direct sum of a unilateral weighted shift operator and a backward weighted
shift operator. Therefore,  the local spectral properties of $S$ can be deduced
from the preceding results.

\section*{Acknowledgments}

Part of this material was presented by the author in the
Mathematics seminar of the Abdus Salam ICTP (August, 2001), in the
Functional Analysis seminar of Rabat (November, 2001), and in the
3rd Linear Algebra Workshop, Bled-Slovenia (June, 2002). The
author expresses his gratitude to the organizers of these
seminars, Professor A. Kuku (Abdus Salam ICTP), Professor O.
El-Fallah (University of Mohammed V, Rabat), and Professors R.
Drnovek, T. Koir, and M. Omlad\`e (3rd LAW, 2002). He acknowledges
Professors P. Rosenthal, and J. Zemanek who were interested in his
talk at Bled. He also thanks Professor L. R. Williams for
providing him reference \cite{wil}. He also thanks Professor M.
Neumann for giving him some remarks and for providing him
some references including \cite{amn}, and \cite{mmn}.
Finally, the author wishes to
express his thanks to Professor P. Aiena for reading and studying
carefully the first version of this paper.


\begin{thebibliography}{99}



\bibitem{aiena} P. Aiena, M. L. Colasante and  M. Gonz\'alez, \newblock{\it Operators which have a closed quasi-nilpotent part},
Proc. Amer. Math. Soc. 130(2002), no. 9, 2701--2710.


\bibitem{amn} P. Aiena, T. L. Miller and  M. M. Neumann, \newblock{\it On a localized single-valued extension property},
Proc. Royal Irish Acad. (to appear).



\bibitem{AA} A. Atzmon, \newblock{\it Power regular operators}, Trans. AMS. V 347 N 8 (1995) 3101-3109.

\bibitem{atzmon} A. Atzmon and M. Sodin, \newblock{\it Completely indecomposable operators
and uniqueness theorem of Cartwright-Levinson type}, J. Functional Anal. 169(1999) 164-188.

\bibitem{BS} C. A. Berger and J. G. Stampfli, \newblock{\it Mapping theorems for the numerical range}, Amer. J.  Math. 89
(1967) 1047-1055.

\bibitem{phd} A. Bourhim, \newblock{\it Points d'\'evaluations born\'ees des op\'erateurs cycliques et
comportement radial des fonctions dans la classe de Nevanlinna}, Th\`ese de Doctorat, Univ. Mohammed V,
Facult\'e des Sciences, Rabat-Maroc, N d'ordre 1924, (2001).


\bibitem{abdou} A. Bourhim, \newblock{\it On the largest analytic set for cyclic operators},
(to appear).

\bibitem{bourhim} A. Bourhim, \newblock{\it On the local spectral properties of weighted shift operators}, The Abdus Salam ICTP Preprint, (2002).


\bibitem{cf} I. Colojoara and C. Foias,
\newblock{\it Theory of generalized spectral operators},
Gordon and Breach, New York 1968.


\bibitem{ln} K. B. Laursen and M. M. Neumann,
\newblock{\it An introduction to local spectral theory}, London Mathematical
Society Monograph New Series 20 (2000).

\bibitem{MC} P. McGuire, \newblock{\it A local functional calculus},
Integral Equations and Operator Theory, 9(1986) 218-236.

\bibitem{mm} T. L. Miller and V. G. Miller, \newblock{\it An operator
satisfying Dunford's condition (C) without Bishop's property
($\beta$)}, Glasgow Math. J. 40(1998) 427-430.

\bibitem{mmn} T. L. Miller, V. G. Miller and M. M. Neumann, \newblock{\it Local spectral properties of weighted shifts}, J. Operator Theory (to appear).



\bibitem{rid}W. C. Ridge, \newblock{\it Approximate point spectrum of a weighted shift}, Trans. Amer. Math. Soc. 147(1970) 349-356.

\bibitem{shi} A. L. Shields, \newblock{\it Weighted shift operators and
analytic function theory}, in Topics in Operator Theory,
Mathematical Surveys, N0 13 (ed. C. Pearcy),  pp. 49-128. American
Mathematical Society, Providence, Rhode Island 1974.

\bibitem{tin} T. Tam, \newblock{\it On a conjecture of Ridge}, Proc. Amer. Math. Soc. 125 (1997), no. 12, 3581--3592.

\bibitem{praha} P. Vrbov\'a, \newblock{\it On local spectral properties of operators in Banach spaces}, Czechoslovak Math. J. 23(1973),
483--492.

\bibitem{wil} L. R. Williams, \newblock{\it Local functional calculus and related results on the single-valued extension
property}, Integral Equations and Operator Theory, accepted for
publication.

\bibitem{yan} L. Yang, \newblock{\it Hyponormal and subdecomposable operators}, J. Functional Anal. 112(1993) 204-217.

\end{thebibliography}
\end{document}